\documentclass[12pt]{amsart}
\usepackage{amsfonts,amssymb}
\usepackage{graphicx}
\usepackage{amscd}
\usepackage{eucal}
\input epsf
\newtheorem{thm}{\bf Theorem}[section]
\newtheorem{cor}[thm]{\bf Corollary}
\newtheorem{lem}[thm]{\bf Lemma}

\theoremstyle{definition}

\theoremstyle{remark}

\def\R{\mathbb{R}}

\def\Z{\mathbb{Z}}

\def\SL{{\rm SL}}
\def\PSL{{\rm PSL}}
\def\<{\langle }
\def\>{\rangle }

\begin{document}
\setlength{\baselineskip}{18pt}
\title[ ]
{ {Symmetries  of  quadratic forms classes and  of quadratic surds
continued fractions.
\\ Part II: Classification of the periods' palindromes} }

\author[F.~Aicardi]{Francesca Aicardi}
\address[F.~Aicardi]{SISTIANA 56 PR (Trieste), Italy}
\email{aicardi@sissa.it }
%
%
%
%

\begin{abstract}
The continue fractions of  quadratic surds are periodic, according
to a theorem by Lagrange.  Their periods may have  differing types
of symmetries. This work relates these types of symmetries to the
symmetries of the classes of the corresponding indefinite
quadratic forms. This allows to classify the periods of quadratic
surds and at the same time  to find, for an arbitrary indefinite
quadratic form, the symmetry type of its class and the number of
integer points, for that class, contained in each domain of the
Poincar\'e model of the de Sitter world, introduced in Part I
\cite{Aic}. Moreover,  we obtain the same information for every
class of forms representing zero, by the finite continue fraction
related to  a special representative of that class. We will see
finally the relation between  the reduction procedure  for
indefinite quadratic forms, defined by the continued fractions,
and the classical reduction theory, which  acquires  a geometrical
description by the results of Part I.
\end{abstract}
\maketitle

\section{Definition of the palindromes}

{\bf Definition 1.} A finite continued fraction
$[\alpha_0,\alpha_1,\dots,\alpha_N]$ is said {\it
palindromic}\footnote{The word palindrome means exactly "that
reads the same backwards and forwards, e.g. "{\sc radar}", or
entire phrases like the Latin riddle "{\sc in girum imus nocte et
consumimur igni}"  ("we go around in the night and we are burnt by
fire"). } iff
\[  \alpha_i=\alpha_{N-i}, \quad \quad i=0,\dots, N. \]

{\bf Definition 2.} A {\it period of length $P$} is a finite
sequence of $P$ natural numbers  that cannot be written as a
sequence of identical  sub-sequences.  For example, [1,2,3,1,2,3]
is not a period.

{\bf Definition 3.} An infinite continued fraction
$[\alpha_0,\alpha_1, \dots]$ is said  {\it periodic} if for some
non negative integer $N$ and some natural $P$ :
\[  \alpha_{N+j}=\alpha_{N+j{\rm mod} P}, \quad \quad \forall j> 0.  \]
It is denoted by:
\[  [\alpha_0, \alpha_1, \dots, \alpha_{N-1},[\alpha_N,\alpha_{N+1},\dots, \alpha_{N+P-1}]],  \]
and  $[a_1,a_2, \dots, a_P]:=[\alpha_N,\alpha_{N+1},\dots,
\alpha_{N+P-1}]$ is its {\it period}, of length $P$.

{\bf Definition 4.} The {\it inverse} of the period of a periodic
continued fraction is obtained by writing the period backwards.
Example: the period $[a,b,c,d]$ is the inverse of the period
$[d,c,b,a]$.

{\bf  Definition 5.}  The period of a periodic continued fraction
is said to be {\it palindromic} if there exists a cyclic
permutation of it such that the permuted period is equal to its
inverse period. For example,  the following periods are
palindromic
\[   [a,a,b,b] \quad , \quad [a,a,b,a,a] \quad   {\rm and} \quad    [a,b,a,c,c]. \]

Note that any period of a sole element is palindromic.

{\bf  Definition 6.}  The period of a periodic continued fraction
is said to be {\it even} if its length is even,  and {\it odd}
otherwise.

{\it  Remark. } If a palindromic period is even, then there exist
at least two different cyclic permutations of it such that the
permuted periods are equal to their inverse.  Example: $[abccba]$
and $[cbaabc]$.

{\bf  Definition 7.}  The period of a periodic continued fraction
is said to be {\it  bipalindromic} if   there is a cyclic
permutation of it such that the permuted period can be subdivided
into two palindromic odd sequences. For example, the periods
\[   \Gamma_1=[a,b,c,b,a,d] \quad   {\rm and} \quad    \Gamma_2=[a,a,a,b] \]
are bipalindromic, since $\Gamma_1$ can be written as
$[(b,c,b)(a,d,a)]$ (or  $[(b,a,d,a,b)(c)]$, etc.)  and $\Gamma_2$
as $[(a)(a,b,a)]$ (or $[(b)(a,a,a)]$).

Note that any period of two different elements is, according to
the definition  above,  bipalindromic.  Hence a non palindromic period
contains at least three different elements.

{\it Remark.} The palindromicity  of a period of $P$ elements
(natural numbers) can be seen as the symmetry, with respect to an
axis, of a plane polygon, whose vertices are labeled by these
natural numbers.  If  $P$ is odd, a vertex must belong to the
symmetry axis, whereas if $P$ even, either no vertices belong to
the symmetry axis (and hence every element has its symmetric
element), or two vertices belong to the symmetry axis (and hence
these two vertices have no symmetric element). This last case
corresponds to bipalindromicity.

\section{The symmetry types of the classes of quadratic forms}

By ${\bf f}$ we denote the triple of integer  coefficients of the
binary quadratic form   $f=mx^2+ny^2+kxy$.

In Part I we called $\mathcal T$ the group, isomorphic to
$\PSL(2,\Z)$, acting on the space of the form coefficients
$(m,n,k)$, and whose action is induced by that of $\SL(2,\Z)$ on
the plane $(x,y)$, where the binary forms are defined.

 The class of the form ${\bf f}=(m,n,k)$ under $\mathcal
T$ is denoted by $C({\bf f})$ or by $C(m,n,k)$.

We recall  the classification of the symmetry types of classes of
indefinite binary quadratic forms (i.e., with discriminant
$\Delta=k^2-4mn<0$),  already introduced in Part I.
  We have considered three commuting involutions, acting in the space of forms 
  and  defining, with ${\bf f}$,  8  forms (see Figure \ref{desit}-III).  In particular, given ${\bf f}=(m,n,k)$:

\begin{enumerate}
\item  the form ${\bf f}_c=(n,m,-k)$ is  the {\it complementary} of the
form ${\bf f}$;

\item   the form $\overline {\bf f}=(m,n,-k)$ is  the {\it conjugate} of the
form ${\bf f}$;

\item   the form ${\bf f}^*=(-n,-m,k)$ is  the {\it adjoint} of the
form ${\bf f}$;

\item the form $\overline{\bf f}^*=(-n,-m,-k)$  is  the {\it antipodal} of the form
${\bf f}$, and is the adjoint of the conjugate (or the conjugate
of the adjoint) of the form ${\bf f}$;

\item  the form $-{\bf f}=(-m,-n,-k)$  is  the {\it opposite} of the form
${\bf f}$, and  it is the complementary of the  adjoint  of the
form ${\bf f}$.
\end{enumerate}

{\it Remarks.}

-- Any form $(m,n,0)$ is {\it selfconjugate}, i.e., $\overline{\bf
f}={\bf f}$.

-- Any form $(m,-m,k)$ is {\it selfadjoint}, i.e.,  ${\bf
f}^*={\bf f}$.

The complementary  of a form ${\bf f}=(m,n,k)$ belongs to the
class of ${\bf f}$, $C(m,n,k)$, while     the  conjugate and/or the adjoint
forms of ${\bf f}$ {\it may} or {\it may not}  belong to  the class of ${\bf f}$.

However (see Proposition 1.2  in Part I) {\sl  if  a class contains a pair of forms related by  some involution, or a
form  which is invariant  under some involution, then the entire class is invariant under that involution}.

Hence, according to the different  symmetries  of forms that  we have considered, there are  exactly five types
of  symmetries of the  classes.

{\bf  Definition 8.}  A class of forms is said to be:

\begin{enumerate}
\item {\it asymmetric} if it contains only pairs of  complementary form;

\item {\it $k$-symmetric} if it contains  only pairs of
complementary forms and conjugate forms or isolated  selfconjugate
forms;

\item {\it $(m+n)$-symmetric} if it contains only pairs of
complementary forms and adjoint forms and isolated selfadjoint
forms;

\item {\it antisymmetric} if it contains only pairs of
complementary forms and antipodal  forms;

\item {\it supersymmetric} if it contains all pairs of
complementary, conjugate, adjoint (and thus antipodal) forms.

\end{enumerate}

{\it Remark.} In term of the symmetries  in the plane $(x,y)$,  given a form $f(x,y)$ ($=f(-x,-y)$) in the class $C(\bf f)$,
\begin{enumerate}
\item  the complementary  form   $f(-y,x)=f(y,-x)$  belongs to  $C(\bf f)$;

\item the  form $f(y,x)$  and its  complementary   form $f(x,-y)=f(-x,y)$ belong to $C(\overline {\bf f})$;

\item the form  $-f(x,y)$ and its complementary  form  $-f(y,-x)=-f(y,-x)$ belong to $C({\bf f}^*)$;

\item the form $-f(y,x)$  and its complementary  form  $-f(x,-y)=-f(-x,y)$ belong to $C(\overline{\bf f }^*)$.
\end{enumerate}

\section{Results}

In this section we enunciate the theorems, which will be proved in
the next sections, using mainly the results of Part I \cite{Aic}.

\subsection{Basic Theorems} Let ${\bf f}:=(m,n,k)$ be a triple of integers such that
$k^2-4mn>0$.

The ordered pair $(\xi^+({\bf f}),\xi^-({\bf f}))$ denotes the
roots (the first with sign plus, the second with sign minus) of
the quadratic equation $ m\xi^2+k\xi+n=0 $:
\begin{equation}\label{roots} \xi^{\pm}({\bf f})=\frac{-k\pm
\sqrt{k^2-4mn}}{2m}.\end{equation}

Suppose $\xi^{\pm}(\bf f)$ be irrational.

\begin{thm}\label{tII-1}  The continued fractions of  roots $\xi^+({\bf f})$ and $\xi^-({\bf f})$
are periodic and their periods, up to cyclic permutations, are one
the inverse of the other\footnote{This theorem was probably
already known to Lagrange, Galois, etc.  A geometrical  proof of
the first part is given in \cite{Ar3}, and of the second part in
\cite{Pav}. }.
\end{thm}

\begin{thm}\label{tII-2} The ordered pair of periods of the continued fractions
of  $(\xi^+({\bf f}),\xi^-({\bf f}))$,
 both periods considered up to cyclic permutations, is an
invariant of the class $C(m,n,k)$. \end{thm}

Because of Theorem \ref{tII-1}, the symmetry properties  of  the
periods of $\xi^+$ and  $\xi^-$ coincide.

{\bf Notation.} In the sequel,
$\Gamma(m,n,k)$ will denote the period of  the continued fraction
of $\xi^+(m,n,k)$.

\begin{thm}\label{tII-2bis}
Any period $s$   defines the class $C(m,n,k)$ such that
$\Gamma(m,n,k)=s$   up to a multiplicative positive integer constant and
up antipodal involution.
\end{thm}

{\bf Definition.}  A form ${\bf f}$ is said to be {\it primitive}
if cannot be written as   $a{\bf f}'$ for another integer form
${\bf f}'$ and  $a>0$.  All forms in the same class are either
primitive or non primitive.  A class of primitive forms is said
primitive.

The  theorems of the next section imply, moreover, the following

\begin{cor}\label{corII-2}

Any period  $s$ of length $P$ defines uniquely the primitive class
$C(m,n,k)$ such that $\Gamma(m,n,k)=s$ if and only if $P$ is odd.
\end{cor}

\subsection{Theorems on  the symmetries of the periods} $\ \ \ $

{\bf Theorems 3.5--9} {\it There is an one-to-one correspondence  between the 5   
symmetry types (in Definition 8) of the  classes  $C(m,n,k)$  and the  5  symmetry types of
their  corresponding periods  $\Gamma(m,n,k)$. }

This one-to-one correspondence is  stated separately  for each type of symmetry  in 
the following 5 theorems.

\begin{thm} \label{tII-3}  The period  $\Gamma(m,n,k)$  is palindromic and even
iff the class $C(m,n,k)$ is $(m+n)$-symmetric. \end{thm}

{\sc Example.}  $m=5,\  n=-7,\  k=9$.  $\Gamma=[1,1,2,2]$.   In
the same class we have:
\[\begin{array} {c|c|c||c}
  m & n & k & \Gamma \\  \hline
  5  & -7 &9 & [1,1,2,2]  \\ \hline
  7   &  -7 & -5 & [1,2,2,1] \\ \hline
  7  &  -5 & 9  & [2,2,1,1] \\ \hline
  5  &  -5  & -11 & [2,1,1,2] \\ \hline
\end{array}
\]
Note that the class contains a pair of adjoint forms and two
self-adjoint forms.

\begin{thm} \label{tII-4}  The period $\Gamma(m,n,k)$ is palindromic and odd iff
the class $C(m,n,k)$ is supersymmetric. \end{thm}

{\sc Example.}  $m=1,\  n=-2,\  k=-3$.  $\Gamma=[3,1,1]$.   In the
same class:
\[\begin{array} {c|c|c||c}
  m & n & k & \Gamma \\  \hline
  1   &  -2 & -3 & [3,1,1] \\ \hline
  1  &  -2 & 3  & [1,1,3 ] \\ \hline
  2   &  -2 & -1 & [1,3,1] \\ \hline
  2  &  -1 & 3  & [3,1,1] \\ \hline
  2  &  -1  & -3 & [1,1,3] \\ \hline
  2  & -2 &1 & [1,3,1]  \\ \hline
\end{array}
\]
Note that the orbit contains two selfadjoint forms, which are
conjugate.

\begin{thm} \label{tII-5}   The period  $\Gamma(m,n,k)$  is bipalindromic iff the
class $C(m,n,k)$ is $k$-symmetric. \end{thm}

{\sc Example.}  $m=1,\  n=-2,\  k=-5$.  $\Gamma=[5,2,1,2]$.   In
the same class:
\[\begin{array} {c|c|c||c}
  m & n & k & \Gamma \\  \hline
  1  & -2 & -5 & [5,2,1,2]  \\ \hline
  1   &  -2 & 5 & [2,1,2,5] \\ \hline
  3  &  -2 & -3  & [1,2,5,2] \\ \hline
  3  &  -2  & 3 & [2,5,2,1] \\ \hline
\end{array}
\]
Note that the orbit contains two pairs of conjugate forms.

{\it Remark.}  The square root of a rational number $\sqrt{p/q}$
has continued fraction with period either odd and palindromic  or
bipalindromic, since it is the root of the equation $qx^2-p=0$,
corresponding to a class of forms either $k$-symmetric or
supersymmetric. This answers a question that Arnold posed in
\cite{Ar2}.

\begin{thm} \label{tII-6}   The period  $\Gamma(m,n,k)$  is non palindromic and odd iff the
class $C(m,n,k)$ is antisymmetric. \end{thm}

{\sc Example.}  $m=5,\  n=-3,\  k=-13$. $\Gamma=[2,1,4]$. In the
same class:
\[\begin{array} {c|c|c||c}
  m & n & k & \Gamma \\  \hline
  5   &  -3 & -13 & [2,1,4] \\ \hline
  5  &  -9 & 7  & [1,4,2] \\ \hline
  3   &  -9 & -11 & [4,2,1] \\ \hline
  3  &  -5 & 13  & [2,1,4] \\ \hline
  9  &  -5  & -7 & [1,4,2] \\ \hline
  9  & -3 & 11 & [4,2,1]  \\ \hline
\end{array}
\]

Note that the orbit contains three pairs of  antipodal forms.

\begin{thm} \label{tII-7}    The period  $\Gamma(m,n,k)$  is non palindromic and even iff the
class $C(m,n,k)$ is asymmetric. \end{thm}

{\sc Example.}  $m=5,\  n=-15,\  k=18$. $\Gamma=[1,2,3,4]$. In the
same class:
\[\begin{array} {c|c|c||c}
  m & n & k & \Gamma \\  \hline
  5   &  -15 & 18 & [1,2,3,4] \\ \hline
  8 &  -15 & -12  & [2,3,4,1] \\ \hline
  8   &  -7 & 20 & [3,4,1,2] \\ \hline
  5  &  -7 & -22  & [4,1,2,3] \\ \hline
\end{array}
\]

In \cite{Ar1}, Arnold posed the question whether the roots of all
quadratic equations of type $x^2+kx+n=0$ are palindromic. The
answer is given by the following corollary.

\begin{cor}\label{corII-1} The continued fractions  of the quadratic surds corresponding to a form whose
class represents 1 have period  either odd and palindromic or even
and bipalindromic.
\end{cor}

\subsection{Theorems on the numbers of representatives}

The theorems below complete the results of Part I and refer to
particular domains of the space of forms, which are defined there.

We will see, proving the theorems above, that  the period
$\Gamma(m,n,k)$ is related  to the cycle (of to half of it)
composed by the representatives of the class $C(m,n,k)$ satisfying
$m>0$ and $n<0$ (thus belonging to $H^0$).

Every class $C(m,n,k)$ has the same symmetry of the cycle of its
representatives in $H^0$. However, a cycle with some of the
symmetries of the classes that we have considered, could possess a
priori some higher symmetry, namely that of an n-gone. The
following theorem excludes this possibility.

\begin{thm}\label{nopol}
The cycle in $H^0$ cannot have  symmetries other  than those of
its class.
\end{thm}

  Let $\Gamma(m,n,k)={[a_1,\dots,a_P]}$.

 {\bf Definition.} If $P$ is odd, define
 \[ \Pi(m,n,k):= \Gamma^2= [a_1,a_2,\dots, a_P, a_{P+1}, \dots,a_p ] \]
where $p=2P$ and $a_{P+i}=a_i$, for $i=1,\dots, P$; otherwise,
$\Pi(m,n,k)=\Gamma(m,n,k)$ and $p=P$.

\begin{thm} \label{tII-8}  Let $(m,n,k)$ be any triple of integers such that
$k^2-4mn>0$ is different from  a square number, and
$\Pi(m,n,k)=[a_1,a_2,\dots, a_p]$. Define
\begin{equation}\label{sum} t_{odd}:=\sum_{i \ odd}^{p} a_i, \quad \quad
t_{even}:=\sum_{i \ even}^{p} a_i, \quad \quad t:=\sum_i^{p} a_i.
\end{equation} Class $C(m,n,k)$ has $t$ points in $H^0$ and  in $H^0_R$,
has $t_{odd}$ points in every domain of $G_A$ and $G_{\bar A}$ and
has $t_{even}$ points  in every domain of $G_B$ and $G_{\bar B}$
(or vice versa).

Moreover $t_{odd}=t_{even}=t/2$ if $\Gamma$ is either odd or even
and palindromic, i.e., if the corresponding form is
supersymmetric, antisymmetric, or $(m+n)$-symmetric.
\end{thm}

In Section 4 of Part I we have seen that if $\Delta$ is equal to a
square number, each class has representatives on the boundaries of
the domains of $C_H$. In particular, Theorem 4.14 says that there
are $k$  distinct classes with discriminant equal to $k^2$. These
$k$ classes have a fixed number of representatives inside each
domain. The following theorems deduce the number of points inside
the domains of $C_H$ of every class and  its symmetry type from
the finite continued fraction of a rational number, related to a
representative  of that class.

{\it Remark.}  The last element of a finite continued fraction is
greater than 1.

{\bf Definition.} We call {\it odd}  continued fraction of a
rational number $r>1$ the finite continue fraction
$[a_1,\dots,a_N]$ of $r$, if $N$ is odd, otherwise the continue
faction $[a_1,a_2,\dots,a_N-1,1]$. Similarly, we call {\it even}
continued fraction of a rational number $r>1$ the finite continue
fraction $[a_1,\dots,a_N]$ of $r$, if $N$ is even, otherwise the
continue faction $[a_1,a_2,\dots,a_N-1,1]$.

Note that the odd (even) continued fraction  of
$r=[a_1,\dots,a_N]$, when $N$ is even (resp., odd), still
represents $r$, indeed
\[[a_1,\dots, (a_N-1),1]
=a_1+\frac{1} {\dots + \frac {1} {(a_N-1)+\frac{1}{1}} }
=a_1+\frac{1} {\dots + \frac {1} {a_N}  }=[a_1,\dots, a_N].
\]

\begin{thm} \label{tII-9}
Let  $k>m>0$ and  $[a_1,\dots,a_L]$ be the even continued fraction
of the rational number $k/m$.  Define \begin{equation}\label{eq4}
t_{odd}:=\sum_{i \ odd}^{L} a_i -1, \quad t_{even}:=\sum_{i \
even}^{L} a_i-1, \quad t:=\sum_i^{L} a_i-1. \end{equation}  The
following statements hold:

i)  Class $C(m,0,k)$ in $H^0$ and $H^0_R$ has $t$ points, in every
domain of $G_A$ and $G_{\bar A}$ has $t_{odd}$ points and in every
domain of $G_B$ and $G_{\bar B}$ it has $t_{even}$ points.

ii) Moreover, $t_{odd}=t_{even}=(t-1)/2$ if $C(m,0,k)$ is
$(m+n)$-symmetric.

 \end{thm}

\begin{thm}\label{t10}

Let $\Delta=k^2$ and $0\le m<|k|$.

i) Class $C(m,0,k)$ is supersymmetric iff  $m=0$ or if $k$ is even
and $m=k/2$ (the discriminant is in this case divisible by 4).

ii) $C(m,0,k)$ is $(m+n)$-symmetric iff  the even  continued
fraction of $k/m$  is palindromic.

iii) $C(m,0,k)$ is $k$-symmetric iff the odd continued fraction of
$k/m$ is palindromic.

iv)  Class $C(m,0,k)$ cannot be antisymmetric.

v) $C(m,0,k)$ is asymmetric iff neither the odd nor the even
continued fraction of $k/m$  are palindromic.

\end{thm}

The Appendix contains  examples of  the  theorems above.

The last section \ref{reduc} is devoted to the reduction theory
for indefinite forms, from the geometrical view point of our
model. We will prove, moreover, the following theorem on the
periodic modular fractions (see eqq.(\ref{mod1}) and
(\ref{mod3})):

\begin{thm} \label{t3} Let $(c_1,c_2,\dots, c_t)$ be the period of the modular  fraction of
a quadratic surd. Then \[ \sum_{i=1}^t c_i =3t \] if   the
corresponding   class is supersymmetric, antisymmetric or
$(m+n)$-symmetric.
\end{thm}

\section{Proofs}

\subsection{Fundamental lemmas}

\begin{figure}[h]
\centerline{\epsfbox{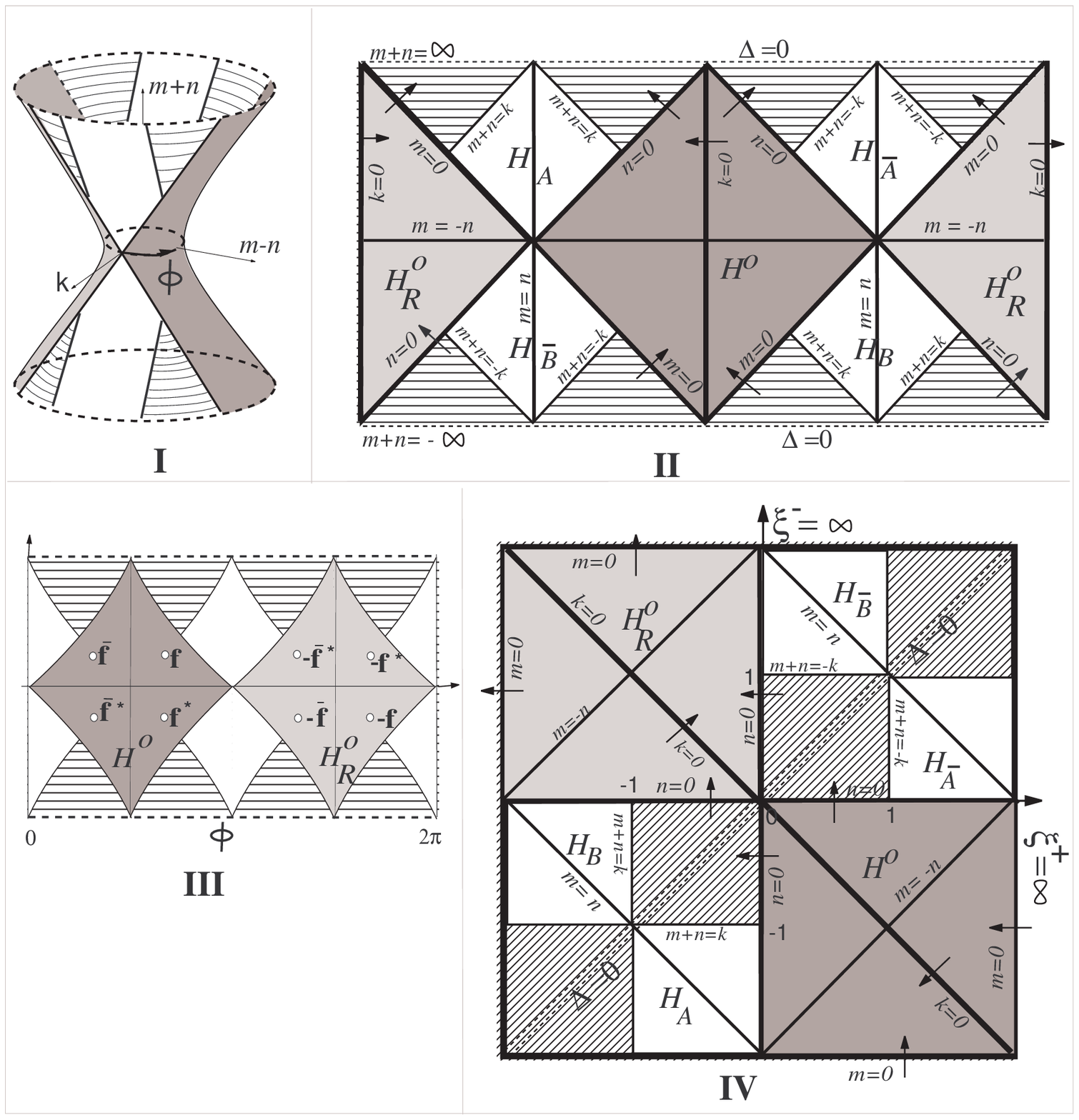}} \caption{}\label{desit}
\end{figure}

In Part I we defined, for  any $\Delta>0$  such that $\Delta \
{\rm mod} \ 4=0,1$:
\[  H_{\Delta}=\{(m,n,k)\in  \R^3 \quad : \quad  k^2-4mn=\Delta \}.\]
This is the space of quadratic forms
\[   f=mx^2+ny^2+kxy \]
with real coefficients and fixed discriminant (see Figure
\ref{desit}-I).

Moreover, we have defined the projection $\mathcal Q$ of the
hyperboloid $H_{\Delta}$   to the open cylinder $C_H$ (Figure
\ref{desit}-III).

\begin{lem}\label{corrisp}  There exists an one-to-one correspondence between
the cylinder  $C_H$   and the space
 \begin{equation} \label{xx} \Xi=\{ (\xi^+({\bf f}),\xi^-({\bf f})) \in \{ \R P^1 \times \R
P^1 \setminus \{ (\xi,\xi)  | \xi \in \R P^1\} ,  \quad {\bf f}\in
H_{\Delta} \}.
\end{equation}
where $\xi^+({\bf f})$ and $\xi^-({\bf f})$ are the roots of the
equation $f=0$ for the variable $\xi=\frac{x}{y}\in \R
P^1$.\end{lem}

{\it Proof.}

 Map $\mathcal Q: H_{\Delta}\rightarrow C_H$ (see
eq. (15) in Part I) is a homeomorphism. Also eq. (\ref{roots})
defines a homeomorphism between $H_{\Delta} $ and $\Xi$. We want,
however, to see explicitly the one-to-one correspondence between
$C_H$ and $\Xi$.

In Figure \ref{desit}-II the cylinder $C_H$ is depicted replacing
the curved segments bounding some of its domains by straight line
segments.

Note that cylinder (\ref{xx}) is obtained by the torus   $\R P^1
\times \R P^1$  minus its diagonal (Figure \ref{desit}-IV). The
circles $c_1$ and $c_2$, frontier  of $C_H$, represent  the limit
points at infinite of the hyperboloid, which coincide with those
of the cone $\Delta=0$. For these limit values of the
coefficients, the roots of the corresponding quadratic equations
tend to a same value. Hence the two circles correspond to the
diagonal $\xi^+=\xi^-$. The root $\xi^+$ vanishes when $n=0$ and
$k>0$, and $\xi^-$ vanishes when $n=0$ and $k<0$.  The values of
$\xi^+$ and $\xi^-=\pm \infty$ are reached when $m=0$, and  they,
too, changes sign when $m$ changes sign. Note that lines $m=0$ and
$n=0$ are the boundaries of domains $H^0$ and $H^0_R$.  The
rhomboidal regions $H^0$ and $H^0_R$ in $C_H$ (Figure
\ref{desit}-III and  Figures 8 and 10 of Part I) are represented
by true rhombi in Figure \ref{desit}-II. These regions are  thus
represented in $\Xi$  by the square regions $\xi^+ \cdot \xi^-<0$,
denoted by $H^0$ and $H^0_R$ in Figure \ref{desit}-IV as well.
Outside $H^0$ and $H^0_R$ coefficients $m$ and $n$ have the same
sign and there are four important domains: $H_A$ and $H_{\bar A}$,
where $m$ and $m$ are positive, $m+n<k$, $k>0$ ($H_A$) and
$m+n<-k$, $k<0$ ($H_{\bar A}$); $H_B$ and $H_{\bar B}$, where $m$
and $n$ are negative, $m+n>k$, $k<0$ ($H_{B}$) and and $m+n>-k$,
$k<0$ ($H_{\bar B}$). They are mapped, respectively,  to the
domains:
\begin{equation}\label{HABt} \begin{array}{llccc}
           H_A & =  \{ (\xi^+,\xi^-): & -1<\xi^+ < 0 &, &  \xi^-<-1 \};  \\
 H_{\bar A} & =  \{(\xi^+,\xi^-): & \xi^+ > 1& , & 0<\xi^-<1  \}; \\
 H_B & =  \{ (\xi^+,\xi^-): & \xi^+ < -1 &, & -1<\xi^-<0 \}; \\
  H_{\bar B} & =  \{ (\xi^+,\xi^-): & 0<\xi^+ <1& , & \xi^->1 \}.
 \end{array} \end{equation}

\begin{figure}[h]
\centerline{\epsfbox{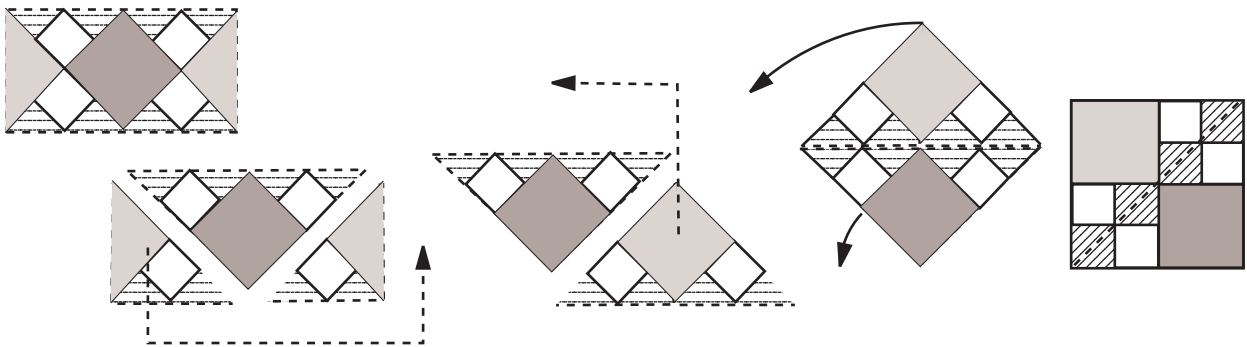}} \caption{Correspondence between
$C_H$ and $\Xi$. }\label{desit2}
\end{figure}

These considerations are sufficient to determine the complete
correspondence. To obtain the `square' of Figure \ref{desit}-IV we
have to cut the `rectangle' of Figure  \ref{desit}-II along the
lines $m=0$, so obtaining  two triangles:  one containing $H^0$
(with circle $c_1$ as base) and the other containing $H^0_R$ (with
circle $c_2$ as base, see Figure \ref{desit2}). Then, we put the
triangle containing $H^0$ above the other, as shown in  figure.
Finally, turn the figure so obtained by $\pi/4$ counterclockwise.
Note that this procedure preserves the continuity. \hfill
$\square$

{\it Remark.}  The complementary form ${\bf f}_c$  of the form
${\bf f}$  in $C_H$ is represented  by a point with the same
ordinate as $\bf f$ and shifted by $\pi$ in the horizontal
direction, whereas the conjugate, the adjoint and the antipodal
forms    are symmetrical with respect to ${\bf f}$ as shown in
Figure \ref{desit}-III.

The following relations hold  among the pairs $(\xi^+,\xi^-)$ of
the triples obtained from the triple ${\bf f }=(m,n,k)$ by all the
considered involutions. \vskip 5 pt
 \centerline{\begin{tabular}{|c|c|c|c|||c|c|c|c| }
\hline
    ${\bf f}$ & $(m,n,k)$ & $\xi^+$ & $\xi^-$  &
       ${\bf f}_c=-{\bf f}^*$ & $(n,m,-k)$ & $-1/\xi^+$ & $-1/\xi^-$
       \\ [ 2
pt]
       \hline
  $\overline  {\bf f}$ & $(m,n,-k)$ & $-\xi^-$ & $-\xi^+$ &
   $\overline{\bf f}_c=-\overline {\bf f}^*$ & $(n,m,k)$ & $1/\xi^-$ &$1/\xi^+$   \\ [ 2
pt]
   \hline
      ${\bf f}^*$ & $(-n,-m,k)$ &$-1/\xi^-$ &$-1/\xi^+$  &
      ${\bf f}_c^*=-{\bf f}$ & $(-m,-n,-k)$ &$\xi^-$ &$\xi^+$    \\ [ 2
pt]
      \hline
 $\overline {\bf f}^*$ &  $(-n,-m,-k)$ &$1/\xi^+$ &$1/\xi^-$  &
 $\overline{\bf f}_c^*=-\overline {\bf f}$ &$(-m,-n,k)$ & $-\xi^+$ &$-\xi^-$    \\  [ 2
pt]   \hline
\end{tabular} }

 \vskip 5 pt

\centerline{Table 1} Therefore, in $\Xi$ the complementary forms
of the forms in  $ H^0$  are obtained moving $H^0$ by a
translation over $H^0_R$, and vice versa, and those outside $H^0$
and outside $H^0_R$ by  moving  the upper-right quarter of $\Xi$
over the lower-left, and vice versa.   The antipodal symmetry,
which is a reflection with respect to the centre of $H_0$ and
$H^0_R$, becomes the reflection with respect to point $(1,-1)$ or
$(-1,1)$, etc.

Remember that we denoted by ${\bf A}=(^{1 \ 1}_{0 \ 1})$, ${\bf
B}=(^{1 \ 0}_{1 \ 1})$, and  ${\bf R}=(^{\ 0 \ 1}_{-1 \ 0})$ the
generators of $\SL(2,\Z)$ acting on the $(x,y)$-plane and by $A$,
$B$, $R$  the corresponding generators of $\mathcal T$.

{\bf Definition.} The operators  $\alpha$, $\beta$ and $\sigma$,
acting on $\R P^1$, and   corresponding to the operators $A,B,R$
of $\mathcal T$, are defined  by
\begin{equation}\label{alphabeta} \alpha (\xi^{\pm}({\bf
f}))=\xi^{\pm}(A({\bf f})), \quad \quad \beta (\xi^{\pm}({\bf
f}))=\xi^{\pm}(B({\bf f})), \quad \quad \sigma(\xi^{\pm}({\bf
f}))=\xi^{\pm}(R({\bf f})) .
\end{equation}

\begin{lem}\label{lemab} The  actions of the operators $\alpha$,
$\beta$ and $\sigma$ on the roots $\xi^{\pm}$  coincide with those
of the inverse of the homographic operators ${\rm A}$, ${\rm B}$
and ${\rm R}$ defined by the generators ${\bf A}$, ${\bf B}$ and
${\bf R}$   of $\SL(2,\Z)$.
\end{lem}

{\it Proof.}  The actions of the  inverse homographic operators
${\rm A}^{-1}$ and  ${\rm B}^{-1}$ (see equations (9) and (12) of
Part I), are:
\begin{equation}  {\rm A}^{-1}: \xi  \rightarrow
\frac{\xi -1}{1}; \ \ {\rm B}^{-1}: \xi  \rightarrow
\frac{1}{-1+1/\xi }; \ \ \ {\rm R}^{-1}={\rm R}: \xi  \rightarrow
- \frac{1}{\xi }.
\end{equation}
On the other hand, if $\xi^{\pm}$ are  the two   roots of
$f(x,y)|_{y=1}=0$, then $\alpha(\xi^{\pm})$ are by definition the
corresponding roots of $f(x+y,y)|_{y=1}=f(x+1,1)=0$, and hence
they are  equal to $\xi^{\pm}-1$.

By definition, $\beta(\xi^+)$ is the first root of the equation
$f(x,x+y)|_{y=1}=0$. Note that
$1/\xi^+=\frac{-k-\sqrt{\Delta}}{2n}$ is the second root, $w^-$,
of the equation $f(1,y)=0$.  Hence  $\beta(w^-)$ is, by the above
definition, the second root of $f(x,y+x)|_{x=1}= f(1,y+1)=0$.  We
thus have $\beta(1/\xi^+)=1/\xi^+-1$. The first root of the equation
$f(x,x+y)|_{y=1}=0$ is thus  equal to
$1/\beta(1/\xi^+)=\frac{1}{-1+1/\xi^+}$. For $\beta(\xi^-)$ the
proof is analogous (exchanging $+$ with $-$, and {\sl first} with
{\sl second}).

Finally, note that $-1/\xi^{\pm}=\frac{k\pm \sqrt{\Delta}}{2n}$
are exactly the first and second roots of $f(-y,x)|_{y=1}$, i.e.,
$-1/\xi^{\pm}=\xi^{\pm}(R({\bf f}))$ and hence they are equal, by
definition, to $\sigma(\xi^{\pm}({\bf f}))$. \hfill $\square$

Denote by $\bar \alpha$ and  $\bar \beta$ the elements
$\alpha^{-1}$ and $\beta^{-1}$.

{\it Remark.}  The  above lemma defines an isomorphism between the
group $\mathcal T$, acting on the space of forms, and the group
generated by $\alpha$, $\beta$ and $\sigma$, acting on the space
of ordered pairs of quadratic equation roots.

The  actions of the
operators $\bar \alpha$ and  $\bar \beta$ are
\[ \bar \alpha (\xi)=\xi+1; \quad \quad \bar \beta(\xi)=\frac {1}{1+1/\xi}.  \]

The following lemma holds for all real numbers,  with  continued
fraction non necessarily periodic and infinite.

 Here
$\alpha^{\rm n}$ and $\beta^{\rm n}$ denote the n-th iterations of
$\alpha$ and $\beta$.

\begin{lem} \label{shift} If $\xi>1$, with continued fraction $\xi=[a,b,c,\dots],$
then
\[  \alpha^{a}(\xi )= [0,b,c,\dots].  \]
If $0<\xi<1$, with continued fraction   $\xi=[0,d,e,g, \dots]$,
then
\[  \beta^{d}(\xi )=[e,g, \dots ]. \]
\end{lem}

 {\it Proof.} By Lemma \ref{lemab},
\[  \alpha(\xi )= \xi -1, \quad {\rm and} \quad \alpha^{a}(\xi )= \xi -a. \] Hence, if
\[ \xi =a+ \frac{1}{b+\frac{1}{c+ \frac{1}{\cdots}}},\]
then
\[ \alpha^a(\xi )= \frac{1}{b+\frac{1}{c+ \frac{1}{\cdots}}}=[0,b,c,\dots]. \]
Moreover, \[  \beta(\xi )=\frac{1}{-1+\frac{1}{\xi }},  \quad {\rm
and} \quad \beta^d(\xi )=\frac{1}{-d+\frac{1}{\xi }}.   \] Hence,
if
\[ \xi =\frac{1}{d+\frac{1}{e+\frac
{1}{g+\frac{1}{\dots}}}} ,\] then
\[ \beta^{d}(\xi )=\frac{1}{-d+ \frac {1} {\frac{1}{d+\frac{1}{e+\frac
{1}{g+\frac{1}{\dots}}}}}} =\frac{1}{-d+
 d+\frac{1}{e+\frac {1}{g+\frac{1}{\dots}}}}=e+\frac
 {1}{g+\frac{1}{\dots}}=[e,g,\dots].\]
\hfill $\square$

\begin{lem} \label{shift2} If $\xi<-1$, with continued fraction
$\xi=-[a,b,c,\dots]$, then
\[ \bar \alpha^{a}(\xi)= -[0,b,c,\dots].  \]
If $-1<\xi<0$, with continued fraction   $\xi=-[0,d,e,g, \dots]$,
then
\[ \bar \beta^{d}(\xi)=-[e,g, \dots ]. \]
\end{lem}

{\it Proof.} It is analogous to the proof of the preceding lemma:
indeed, by Lemma \ref{lemab} and the Remark following it,
\[   \bar \alpha ^a = \quad \alpha^{-a}(\xi )= \xi +a; \]
so, if  $\xi=-[a,b,c,\dots]$, i.e.,
\[ \xi=-a-\frac{1}{b+\frac{1}{c+ \frac{1}{\cdots}}},\ \]
then  \[ \xi+a=-\frac{1}{b+\frac{1}{c+ \frac{1}{\cdots}}}.\ \]
Moreover,
\[  \bar \beta^d(\xi)= \beta^{-d}(\xi)= \frac{1}{d +1/\xi};  \]
so, if  $\xi=-[0,d,e,g, \dots]$,  i.e.,
\[  \xi= -\frac{1}{d+\frac{1}{e+ \frac{1}{g+\frac{1}{\cdots}}}} \ , \]
then
\[  \frac{1}{d +1/\xi}= -\frac{1}{e+ \frac{1}{g+\frac{1}{\cdots}}}. \]
\hfill $\square$

We group in the following Lemma some observations.

Let a black arrow from the point ${\bf f}$ to the point $ {\bf g}$ in
$C_H$  indicate that $A{\bf f}= {\bf g}$, and a white arrow  that
$B{\bf f}= {\bf g}$.

\begin{lem}\label{arrows}  If the points ${\bf f}$ and $ {\bf g}$ in $C_H$
are joined by an arrow from ${\bf f}$ to $ {\bf g}$, then

1)  the points ${\bf f}^*$ and $ {\bf g}^*$ in $C_H$, symmetric
with respect to the horizontal line of points ${\bf f}$ and $ {\bf
g}$, are joined by an arrow from ${\bf g}^*$ to $ {\bf f}^*$ of
the opposite color (see Figure \ref{symm}-I,II);

2) the points $\overline{\bf f}$ and $\overline{\bf g}$ in $C_H$,
symmetric with respect to the the vertical  line $k=0$  of points
${\bf f}$ and $ {\bf g}$, are related by an arrow from $\overline
{\bf g}$ to $\overline{\bf f}$ of the same color (see Figure
\ref{symm}-II,III);

3)  the points $\overline{\bf f}^*$ and $\overline{\bf g}^*$ in
$C_H$, symmetric with respect to the centre of $H^0$   of points
${\bf f}$ and $ {\bf g}$,  are related by an arrow from
$\overline{\bf f}^*$ to $\overline{\bf g}^*$  of the opposite
color (see Figure \ref{symm}-II-IV).

\end{lem}

{\it Proof.}  The proof of  the corresponding identities:
\begin{equation}\label{ident} \begin{array} {l l l}
1)   &
 {(A  {\bf f})}^* = \bar B {\bf f}^*;  &    {(B{\bf f})}^* =
\bar A {\bf f}^* ; \\  2)  &
    \overline {A  {\bf f}} = \bar A \ \overline {\bf f};
    &
 \overline{B{\bf f}} = \bar B \ \overline {\bf f}  ; \\
3)  &
    (\overline {A  {\bf f}})^* = B \overline {\bf f}^*;
    &
 (\overline{B{\bf f}})^* = A \overline {\bf f}^* . \end{array}
 \end{equation}
is given in Part I, Lemma 1.3.

\subsection{Classes non representing zero}

Because of the correspondence stated by Lemma \ref{corrisp}, the
same symbols will denote the regions in $C_H$ and in  $\Xi$, like
in Figure \ref{desit}.

\subsubsection*{Proof of   Theorems  \ref{tII-1} and \ref{tII-2}. }
We prove initially  that in $H^0$   the continued fractions of
$\xi^+$ and $\xi^-$ are periodic and their periods are one the
inverse of the other.

{\bf Definition.}  We call {\it principal point} a point ${\bf f}$
belonging to a cycle in $H^0$ such that its consecutive point is
reached by $A$ iff ${\bf f}$ is reached by $B$ from its preceding
point, or vice versa.

 Let   ${\bf h}$ be a principal point  in $H^0$. It
belongs, by Theorem 4.4 and Corollary 4.12 of Part I, to a cycle
$\gamma_{\bf h}(T_1,\dots T_t)$, being $T_i$ equal either to $A$
or to $B$. Let $T_1=A$.

Hence we can write $T{\bf h}={\bf h}$, where
$T=B^{a_p}A^{a_{p-1}}\dots B^{a_2}A^{a_1}$, grouping together in
$p$ groups,  subsequent operators $T_i$ of type $A$ or of type
$B$, so that $\sum_{i=1}^p a_i=t$. The cycle of ${\bf h}$ will be
thus $\gamma_{\bf h}(A^{a_1},B^{a_2},\dots,A^{a_{p-1}},B^{a_p})$,
indicating that ${\bf h}_{a_1}=A^{a_1}{\bf h}$,  $\quad {\bf
h}_{a_1+a_2}=B^{a_2}A^{a_1}{\bf h}$, etc, till  ${\bf h}_{t}=
T{\bf h}={\bf h}$.  Let $(\xi^+,\xi^-)$ be the pair of roots
associated to ${\bf h}$. By Lemmas \ref{corrisp} and \ref{lemab},
there is an analogous operator to $T$, say $\tau$, obtained from
$T$ by translating $A$ into $\alpha$ and $B$ into $\beta$,
satisfying:
\begin{equation}\label{tau} \tau (\xi^+)=\xi^+, \quad \quad \tau
(\xi^-)=\xi^-.\end{equation} We  observe now that, since ${\bf h}$
is in $H^0$, $\xi^+$ is positive and $\xi^-$ is negative. Write
$\xi^+=[b_1,b_2,\dots]$. Applying $A^{a_1+1}$ to ${\bf h}$, we
exit from $H^0$, by Lemma 4.6 of Part I;  so, we obtain
$\alpha^{a_1}\xi^+=[b_1-a_1,b_2,\dots]=[0,b_2,\dots]$ and hence
$b_1=a_1$. Analogously, applying  $B^{a_2+1}$ to $A^{a_1}{\bf h}$,
we exit from $H^0$; so $\beta^{a_2}\circ
\alpha^{a_1}(\xi^+)=[b_3,b_4,\dots]$, i.e., $b_2=a_2$. In this way
we obtain $b_i=a_i$ for $i=1,\dots,p$. But, at the end of the
cycle, we get, by eq. (\ref{tau})
\[ \tau
(\xi^+)=[b_{p+1},b_{p+2}, \dots]= [a_{1},a_{2}, \dots, a_p,
b_{p+1}, b_{p+2},\dots ]. \] Similarly, applying by recurrence to
$\xi^+$, for every natural $j$, the $j$-th iterate $\tau^j$ of
$\tau$:
\[ \tau^j (\xi^+)=[b_{jp+1},b_{jp+2},\dots] = [a_{1},a_{2}, \dots, a_p, a_{1}, a_{2},\dots ],  \]
we obtain  that $[a_{1},a_2, \dots, a_{p-1},a_p]$ are the first
$p$ elements of the continued fraction obtained  from $\xi^+$
canceling out the first $jp$ elements,   so concluding
\[  \xi^+=[[a_1,a_2,a_3,\dots, a_p]]. \]
Note that $p$ is even, since $T$ begins with a power of $B$ and
ends with a power of  $A$.

We consider now $\xi^-$. In this case it is convenient to write
the right equation in (\ref{tau}) as
\[  \xi^-= \tau^{-1}(\xi^-), \]
where $\tau^{-1}= \bar \alpha^{a_1} \circ \bar \beta^{a_2}\cdots
\bar \alpha ^{a_{p-1}} \circ \bar \beta^{a_p}$. Again by Lemma 4.6
of Part I, we know that, applying  $\bar A$ to ${\bf h}$, we exit
from $H^0$, i.e,  $\bar \alpha (\xi^-) = \xi^- +1>0$, and hence
$0>\xi^->-1$, i.e., $\xi^-= -[0,c_1,c_2,\dots]$. For analogous
arguments, we will find, using  again Lemma \ref{shift2}, that
$\bar \beta^{a_p} (\xi^-)=-[c_2,c_3, \dots]$ and hence $c_1=a_p$.
Analogously, applying  $\bar A^{a_{p-1}}$ to $\bar B^{a_1}{\bf
h}$, we obtain  $\bar \alpha^{a_{p-1}}\circ \bar
\beta^{a_p}(\xi^-)=-[c_3,c_4,\dots]$, i.e., $c_2=a_{p-1}$. In this
way, we obtain $c_i=a_{p+1-i}$, for $i=1,\dots,p$. So, at the end
of the cycle, we get, by  eq. (\ref{tau})
\[\tau^{-1} (\xi^-)=-[c_{p+1},c_{p+2}, \dots]= -[a_{p},a_{p-1}, \dots, a_2,a_1, c_{p+1}, c_{p+2},\dots ], \]
and hence also that $[a_{p},a_{p-1}, \dots, a_2,a_1,]$ are the
first $p$ elements of the continued fraction of any iterate of
$\tau^{-1}$ on $\xi^-$:
\[ \tau^{-j} (\xi^-)=-[c_{jp+1},c_{jp+2},\dots] = -[a_{p},a_{p-1}, \dots, a_1, a_{p}, a_{p-1},\dots ],  \]
so concluding
\[  \xi^-= -[0,[a_p,a_{p-1},\dots,a_2, a_1]]. \]

We have proved that the roots $\xi^\pm({\bf h})$  are periodic and
their periods are one the inverse of the other when ${\bf h}\in
H^0$ is a {\it principal point}, i.e., the operator $T$ of
$\mathcal T^+$ satisfying $T {\bf h}={\bf h}$, starts by  $A$ and
ends by $B$ (or vice versa).

It is now clear that a non principal point  in  $ H^0$, in the
cycle of ${\bf h}$ between ${\bf h}$ and ${\bf h}_{a_1}$ (whenever
$a_1>1$), is defined by
\[  {\bf h}_{j}=  A^j {\bf h},  \]
for some $j<a_1$, and satisfies
\[ \xi^+({\bf h}_j)=[a_1-j,a_2,\dots]=[a_1-j,[a_2,a_3,\dots, a_p,a_1]],
\]
i.e., it has the same period as $\xi^+({\bf h})$, being the period
defined up to cyclic permutations. Moreover
\[ \xi^-({\bf h}_j)=-[j,a_p,a_{p-1},\dots]=-[j,[a_{p}, a_{p-1},\dots, a_1]].
\]
Since we may repeat the reasoning for any point of the cycle
between two principal points,  the periods of the continued
fractions of  $\xi^\pm({\bf h})$ are one the inverse of the other
for every ${\bf h}\in H^0$.

To complete the proof of  Theorem \ref{tII-1},  we have to
consider points outside $H^0$.  Every point belongs to an orbit
and every orbit possess, by the results of Part I, a
representative inside $H^0$. Any point of the orbit can be thus
written as ${\bf p}=T{\bf h}$, with ${\bf h}\in H^0$ and $T\in
\mathcal T$. Now, every element $T$ of the group can be written as
a finite product of the generators $A$, $B$ and $R$.  By Lemma
\ref{lemab}, we translate  $T$ into $\tau$, composed by the
corresponding generators $\alpha$, $\beta$ and $\sigma$.  The
action of each of these generators on a continued fraction affects
evidently only its initial elements, and hence $\tau$ affects only
a finite initial part of $\xi^{\pm}({\bf h})$, so that the periods
of $\xi^\pm(\tau({\bf h}))$  remain unchanged, up to cyclic
permutations.

For instance,  the periods of the continued fractions of
$\xi^\pm({\bf h_c})$ are the same as those of $\xi^\pm({\bf h})$,
because of the relations shown in  Table 1, since ${\bf h}_c=R
{\bf h}$.    We have completed the proof of Theorem  \ref{tII-1}
and at the same time that of Theorem \ref{tII-2} \hfill $\square$

{\it Remark.}  The above proof shows, for every class $C(m,n,k)$
of indefinite forms with discriminant different from a square
number, a correspondence between the cycle in $H^0$ (Theorem 4.13
of Part I) and the periods of continued fractions of
$\xi^{\pm}(m,n,k)$. As we remarked,  the number $p$ of groups of
operators of type $A$ and $B$ which alternate in the cycle is
necessarily even.  The length $P$ of the period of the continued
fraction can be, however, odd. In this case the cycle corresponds
indeed to the double of the period, and $p=2P$, as we will see in
detail in Theorems \ref{tII-4} and \ref{tII-6}.

\subsubsection*{Proof of Theorem \ref{tII-2bis}} By the proof of
Theorem \ref{tII-1}, a root $\xi^+({\bf f})$ or  $\xi^-({\bf f})$
is immediately periodic (i.e., there are no elements different
from 0 before the period) only if ${\bf f}$ is a principal point
of $H^0$ or of $H^0_R$, and every class of forms, with
discriminant different from a square number, has a principal
representative in $H^0$ and in $H^0_R$.

Note that $\xi^+$ is positive and $\xi^-$ is negative in $H^0$,
and vice versa in $H^0_R$.  Remember also that  inverting the
order of the pair of roots (and hence of periods) corresponds to
inverting the signs of the coefficients of the equation (i.e., of
${\bf f}$).

Given a sequence $s=[a_1,\dots,a_P]$, we firstly suppose that
$[[a_1,\dots,a_P]]$ is  the first root of a quadratic equation
with integer coefficients. To know such equation, we write:
\[  \xi= a_1+\frac{1}{a_2+\frac {1}{\cdots +\frac {1}{a_P+1/\xi}}}.\]
Call $(m,n,k)$ its coefficients and ${\bf f}=(m,n,k)$.

If we suppose that the first root is less than 1, i.e., it is
equal to $[0,[a_1,a_2,\dots,a_P]]$, we obtain another equation:
\[  \xi= \frac {1}{a_1+\frac{1}{a_2+\frac {1}{\cdots +\frac {1}{a_P+\xi}}}},  \]
which will correspond to  $\overline {\bf f}^*$ (see Table I).

Alternatively, supposing the first root be negative and greater
than $-1$, i.e., equal to $-[0,[a_1,a_2,\dots,a_P]]$, we obtain
another equation which correspond to $-{\bf f}^*$.

Finally, supposing the first root be negative and less than $-1$,
i.e., equal to $-[[a_1,a_2,\dots,a_P]]$, we obtain another
equation which correspond to $-\overline {\bf f}$.

Since $-{\bf f}^*={\bf f}_c$, and $-\overline {\bf
f}=\overline{\bf f}^*_c$, the four triples of coefficients that we
have obtained belong to only two classes, related by the antipodal
symmetry. \hfill $\square$

\begin{figure}[h]
\centerline{\epsfbox{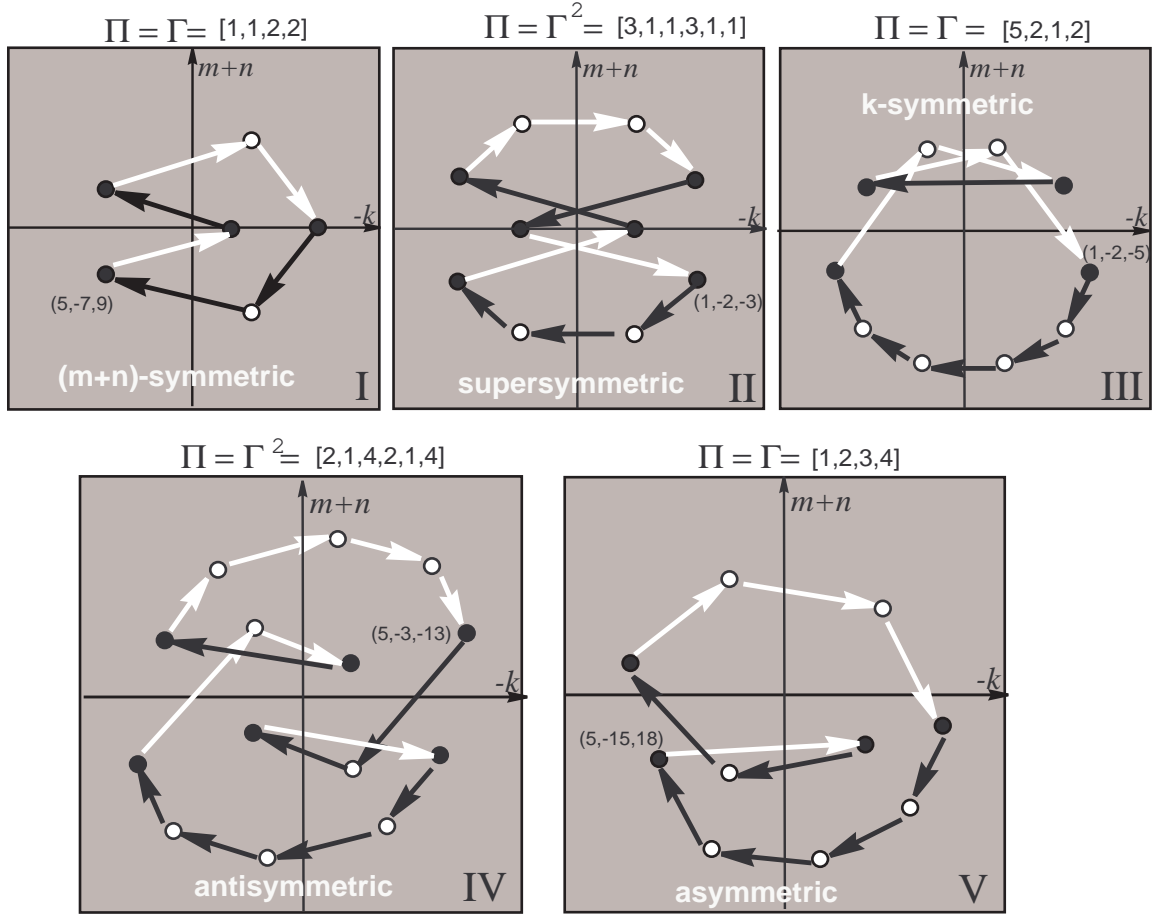}} \caption{The elements of the
periods are equal to the number of arrows between two principal
points (black dots).}\label{symm}
\end{figure}

We give here the proofs of the theorems on the symmetries of the
periods.  We illustrate   theorems \ref{tII-3}-\ref{tII-7}
  in Figure \ref{symm}. The cycles in
$H^0$ correspond indeed  to the periods related to the forms
considered in the  examples.  In these figures the small circles
indicates the forms.  Black circles correspond to principal
points. The black arrow indicate the operator $A$ and the white
arrow the operator $B$.

\begin{lem}\label{2points} If a cycle contains two points related by some  symmetry,
the cycle possesses that symmetry.
\end{lem}

{\it Proof.} By Lemma 4.6 of Part I, a point of a cycle determines
uniquely both the successive  points and the preceding one, and,
hence,  all points of the cycle. Because of relations stated by
Lemma \ref{arrows},  two symmetric points have symmetric arrows
entering and exiting from them, and hence two symmetric points
determine the symmetries of their neighboring in the cycle, and so
the symmetry of the entire cycle. \hfill $\square$

Note that the symmetries of a cycle in $H^0$ concerns {\it only
its points} and the segments joining them: the directions of the
arrows and their colors are necessarily related by the rules given
by Lemma \ref{arrows}.

\begin{lem}\label{vert} There are no vertical arrows in a cycle which contains more than two points.
\end{lem}
{\it Proof.}  Suppose a  vertical white arrow  relates  two
symmetric points ${\bf f}$ and ${\bf f}^*$.  By the lemma above, a
black arrow must relate ${\bf f}^*$ to $({\bf f}^*)^*={\bf f}$. So
either these points form a cycle, or they cannot be joined by an
arrow. \hfill $\square$

\begin{lem}\label{selfadj}
A cycle $(m+n)$-symmetric contains exactly two selfadjoint points.
The selfadjoint points are principal.
\end{lem}
{\it Proof.} Let ${\bf f}$ and ${\bf f}^*$  be two different
points of the cycle (so they lie in the same vertical line).  Let
$ {\bf g}$ be the successive of ${\bf f}$  (i.e., either $A{\bf
f}= {\bf g}$ or $B{\bf f} = {\bf g}$). The symmetric ${\bf g}^*$
is in the cycle and is related to ${\bf f}^*$ by $\bar B{\bf f}^*$
or $\bar A{\bf f}^*$. Continuing with the successive points
according to the arrows after $ {\bf g}$, and their symmetric,
according to the inverse arrows after $ {\bf g}^*$, we have to
close that part of the cycle, from ${\bf f}$ to  ${\bf f}^*$. But
since there are no arrows between adjoint points by Lemma
\ref{vert}, we must have, for some pair of adjoint points $ {\bf
h}$ and $ {\bf h}^*$, that $A {\bf h}={\bf j}$ and $\bar B {\bf
h}^*={\bf j}$, with $ {\bf j}$ selfadjoint (belonging to the line
$(m+n)=0$). Point ${\bf j}$ is thus a  principal point. In the
remaining part of the cycle, from ${\bf f}^*$ to  ${\bf f}$, there
is a selfadjoint point for the same argument. We have now to prove
that there  are no other selfadjoint points.  By the above
argument, starting by any pair of adjoint points an reaching the
successive pairs following the arrows in both directions,  when we
arrive the second selfadjoint point, we close  the  cycle. Since a
point of the cycle cannot be visited twice by results of Part I,
no other selfadjoint points are possible. \hfill $\square$

\begin{lem}\label{1point}
If an orbit contains a selfadjoint point,  is symmetric with
respect to the axis $(m+n)=0$.
\end{lem}

{\it Proof.} Let ${\bf h}={\bf h}^*$ be the selfadjoint point. For
the arguments of the preceding Lemma, if the successive of ${\bf
h}$ is $A {\bf h}$, then the preceding of ${\bf h}$ is $\bar B
{\bf h}$. $A{\bf h}$ and $\bar B{\bf h}$ are   a pair of adjoint
points.  By Lemma \ref{2points}, the orbit is at least
$(m+n)$-symmetric. \hfill $\square$

 \subsubsection*{ Proof of Theorem \ref{tII-3}} (See Figure  \ref{symm}-I).

Suppose the class  be $(m+n)$-symmetric. This means that the cycle
is invariant under reflection with respect to the horizontal axis
($m+n=0$). Choose one of the two  selfadjoint points, $ {\bf h}$,
of the cycle (it exists by Lemma \ref{selfadj}).  Consider all
points of the cycle between $ {\bf h}$ and  the second selfadjoint
point, $ {\bf h}'$, following the arrows. This procedure defines
an operator $M\in \mathcal T^+$ such that
\[ M {\bf h}= {\bf h}'.\]  The points between  $ {\bf h}'$ and $ {\bf h}$, forming the
other part of the cycle, are the adjoint of the points from $ {\bf
h}$ to $ {\bf h}'$. By Lemma \ref{arrows},
\[   \hat M {\bf h}= {\bf h}',  \] where $\hat M$ is obtained by $M$
exchanging $A$ with  $\bar B$, and $B$ with $\bar A$.
Equivalently, we may write
\[   \hat M^{-1} {\bf h}'= {\bf h}, \]
where  $\hat M^{-1}$  is obtained by $M$ inverting the order of
the factors and exchanging $A$ with $B$.  But thus  $\hat M^{-1}
=M^\vee$, i.e., the transpose of $M$, since the transpose of $A$
and $B$ are $B$ and $A$, respectively.  Hence
\[ M M^\vee {\bf h}= {\bf h}.   \]
So, if $M$ is  the product of $q$  groups of generators of type
$A$ and $B$: $M=A^{a_1}B^{a_2}\cdots A^{p-1}B^{q}$, the operator
$T=M M^\vee$, defining the cycle, is he product
$A^{a_1}B^{a_2}\cdots A^{q-1}B^{q}A^{q}B^{q-1}\cdots
B^{a_2}A^{a_1}$. The roots $\xi^\pm( {\bf h})$, satisfying: $\tau
(\xi^{\pm}) = \xi^{\pm}$ are thus
\[  \xi^+=[[a_1,a_2,\dots, a_q,a_q,\dots a_2,a_1]], \quad \xi^-=[0,[a_q,a_{q-1},\dots, a_1,a_1,a_2, \dots, a_q]].  \]
Their periods have length $P=2q$ and are palindromic.

On the other hand, given a  period  $\Gamma$ of length $P$, we
associate to it a cycle, as we explained in the proof of Theorem
1, using Lemma \ref{shift}.  Now, if  $\Gamma$ is even and
palindromic, also the sequence of groups  of operators  $A$ and
$B$ is even and palindromic in $T$, satisfying $T {\bf h}= {\bf
h}$, for some principal point $ {\bf h}$. We thus can write $T= M
M^\vee $, and
 \[{\bf h}= M M ^\vee  {\bf h}= (M^\vee)^{-1} M^{-1}{\bf h}. \]
  Since ${M^\vee}^{-1}=\hat M$, and $\hat M^{-1}= M^\vee$, we
  have, by the second equality above
\[  {\bf h}=\hat M M^{-1}{\bf h}, \quad  \quad {\rm and }\quad \quad   {\bf h}^*= M  M^\vee {\bf h}^*. \]
Hence $ {\bf h}= {\bf h}^*$ is selfadjoint.   By Lemma
\ref{1point}, the cycle is $(m+n)$-symmetric. \hfill $\square$

{\it Remark.} For the forms of such a class, $\Pi=\Gamma$ and
$p=P=2q$.

\subsubsection*{Proof of Theorem \ref{tII-4}} (See Figure
\ref{symm}-II).

 Suppose the class be supersymmetric. This means that the cycle is invariant
under reflection with respect to the horizontal axis ($m+n=0$) as
well as with respect to the vertical axis ($k=0$). Since a
supersymmetric orbit has, in particular, the symmetry of a
$(m+n)$-symmetric orbit, it has two selfadjoint points, by Lemma
\ref{2points}.  Since these points lie on the horizontal line
$(m+n)=0$, and the cycle is supersymmetric, these points must be
symmetric with respect  to the vertical axis $k=0$, i.e., they are
conjugate. Call them ${\bf h}$ and $\overline {\bf h}$.  Since
${\bf h} ^*={\bf h}$ and $\overline{\bf h}^*= \overline {\bf h}$,
being these points selfadjoint, they must be principal points. Let
$A {\bf h}$ and $\bar B {\bf h}$ be the successive and preceding
point of ${\bf h}$. The successive and preceding of $\overline
{\bf h}$ are, by Lemma \ref{arrows}, $B \overline {\bf h}$ and
$\bar A \overline {\bf h}$, and
\[    \overline{ A {\bf h}}= \bar A \overline {\bf h}, \quad \quad  \overline{ \bar B {\bf h}}=B \overline {\bf h}.  \]
If we consider the successive points of the cycle, after ${\bf h}$
and $A {\bf h}$, we must meet the point $\overline {\bf h}$, and
necessarily, before it, $\bar A \overline {\bf h}$.   To every
point ${\bf f}= T {\bf h}$, successive to $ {\bf h}$, there
correspond, in the same arc of cycle from ${\bf h}$ to $\overline
{\bf h}$, a symmetric point, $\overline{\bf f}= \bar T \overline
{\bf h}$, where $\bar T$ is obtained by $T$ substituting each $A$
by $\bar A$, and each $B$ by $\bar B$. Therefore  we can write
\[ \overline {\bf h}= M {\bf h}, \quad \quad  {\bf h}= \bar M \overline {\bf h},  \]
obtaining $M=\bar M^{-1}$. I.e., $M$ consists in  a palindromic
sequence
\[  M=A^{a_1}B^{a_2}\cdots  B^{a_{p-1}} A^{a_1}.\]
Note that the sequence have to end by $A$ (or by $B$) if it starts
by $A$ (respectively, by $B$).  Hence the number of groups of $A$
and $B$ in  $M$ is odd.

Using Theorem \ref{tII-3}, to close the cycle from $\overline {\bf
h}$ to ${\bf h}$, we multiply $M$ by $M^\vee$, which is composed
by the same sequence of $q$ groups of generators, in inverse order
and exchanging $B$ with $A$.  But since $M$ is itself palindromic,
it is  invariant under the inversion of the order of factors, and
hence $M^\vee$ is obtained by $M$ simply exchanging $A$ with $B$.

Using Lemmas \ref{corrisp} and \ref{shift} as in  proof of Theorem
\ref{tII-1}, we thus obtain  that the continued fractions of
$\xi^\pm ( {\bf h})$  are odd, 2 periods of them corresponding in
fact to a cycle in $H^0$.

Conversely, if we have a continued fraction with odd palindromic
period \[[a_1, a_2, \dots ,a_{q-1},a_q,a_{q-1},\dots  a_2,a_1]\]
we associate to it an operator $M=A^{a_1}B^{a_2}\dots
B^{a_2}A^{a1}$, and  solving $M M^\vee {\bf h}={\bf h}$, we find a
selfadjoint form, and, by the symmetry of $M$, the cycle of ${\bf
h}$ results to be supersymmetric. \hfill $\square$

{\it Remark.} For the  forms of such a class we thus have
$\Pi=\Gamma^2$ and $p=2P=2(2q+1)$.

\subsubsection*{Proof of Theorem \ref{tII-5}} (See Figure
\ref{symm}-III).

 Suppose the
class be $k$-symmetric. This means that the cycle is invariant
under reflection   with respect to the vertical axis ($k=0$).
Choose any pair in the cycle of conjugate points, and call them
${\bf h}$ and $\overline {\bf h}$.   Let $A {\bf h}$ and $\bar B
{\bf h}$ be the successive and preceding point of ${\bf h}$. The
successive and preceding of $\overline {\bf h}$ are, by Lemma
\ref{arrows}, $B \overline {\bf h}$ and $\bar A \overline {\bf
h}$, and
\[    \overline{ A {\bf h}}= \bar A \overline {\bf h}, \quad \quad  \overline{ \bar B {\bf h}}=B \overline {\bf h}.  \]
If we consider the successive points of the cycle, after ${\bf h}$
and  $A {\bf h}$, we must meet the point $\overline {\bf h}$, and
necessarily, before it, $\bar A \overline {\bf h}$.   To every
point ${\bf f}= T {\bf h}$, successive to $ {\bf h}$, there
correspond, in the same arc of cycle from ${\bf h}$ to $\overline
{\bf h}$, a symmetric point, $\overline{\bf f}= \bar T \overline
{\bf h}$, where $\bar T$ is obtained by $T$ replacing each $A$ by
$\bar A$, and each $B$ by $\bar B$. We write finally
\[ \overline {\bf h}= M {\bf h}, \quad \quad  {\bf h}= \bar M \overline {\bf h},  \]
obtaining $M=\bar M^{-1}$. I.e., $M$ consists in a palindromic
sequence
\[  M=A^{a_1}B^{a_2}A^{a_3}\cdots B^{a_{q-1}} A^{a_q}B^{a_{q-1}}\cdots  A^{a_3} B^{a_2} A^{a_1}.\]
Note that the sequence have to end by $A$ (or by $B$) if it starts
by $A$ (respectively, by $B$).  The central group is $A^{a_q}$ if
$p$ is odd, $B^{a_q}$ is $q$ is even.  The number of groups of $A$
and $B$ in $M$ is $2q+1$, and hence is odd. Applying the same
argument to the arc of circle from $\overline {\bf h} $ to ${\bf
h}$, we obtain  ${\bf h} = N \overline  {\bf h}$, where $N$ is (if
$M$ starts with $A$),
 \[  N=B^{b_1}A^{b_2}B^{b_3}\cdots  B^{b_{r-1}} A^{b_r}B^{b_{r-1}}\cdots B^{b_3} A^{b_2} B^{b_1},\]
and, as before, the central group is $A^{b_r}$ is $r$ is even,
$B^{b_r}$ is $r$ is odd.

We thus obtain
\[  {\bf h} =  NM {\bf h},  \]
where $T=NM$, defining the cycle, is  composed  by two sequences
palindromic and odd:
\[   [(a_1, \dots, a_q, \dots a_1)(b_1, \dots, b_r, \dots , b_1)]. \]
The resulting sequence is, by Theorem \ref{tII-1}, the period of
the root of the quadratic equation associated to ${\bf h}$, and
this period is evidently even and bipalindromic.

Conversely, having a period bipalindromic, we   subdivide it into
two palindromic odd periods, of length $2q+1$ and $2r+1$,
$[a_1,\dots,a_q,\dots a_1]$ and $[b_1, \dots, b_r,\dots, b_1]$ and
build, for instance, the operators $M=A^{a_1}B^{a_2}\dots A^{a_1}$
and $N=B^{b_1}A^{b_2}\dots B^{b_1}$. By Theorem \ref{tII-1}, this
is the period of the  root $\xi^+({\bf f})$, where ${\bf f}$
satisfy:
    \[ {\bf f} =  NM{\bf f} \]
Using Lemma \ref{arrows} we immediately see that  points ${\bf f}$
and $M{\bf f}$, which are in the cycle, are conjugate, and all
points $A{\bf f}, A^2{\bf f}$, $\dots$, $A^{a^1}{\bf f}$,
$A^{a_1}B{\bf f}$, etc, belonging to that part of the  cycle,
defined by $M$, are conjugate of the corresponding points of the
same part of the cycle $\bar A M{\bf f},\bar A^2 M{\bf f}$,
$\dots$, $\bar A^{a^1} M{\bf f}$, $\bar A^{a_1}\bar B M{\bf f}$,
etc. Note that, if $a_q$ is even, there is a central point, on the
axis $k=0$, which is self-conjugate. Similarly the pairs of points
obtained by the second part of the cycle, from $M{\bf f}$ to $N
M{\bf f}={\bf f}$, are conjugate.  So we obtain that the cycle,
being symmetric with respect to the vertical axis $(k=0)$, is
$k$-symmetric. \hfill $\square$

{\it Remark.} In this case  $\Pi=\Gamma$ and $p=P=2(r+q)+2$.

\subsubsection*{Proof of Theorem \ref{tII-6}.} (See Figure  \ref{symm}-IV)

 Suppose the class be antisymmetric. This means that the cycle is invariant under
reflection   with respect to the centre of $H^0$. Choose any pair
in the cycle of antipodal points, and call them ${\bf h}$ and
$\overline {\bf h}^*$.   Let $A {\bf h}$ and $\bar B {\bf h}$ be
the successive and preceding point of ${\bf h}$. The successive
and preceding of $\overline {\bf h}^*$ are, by Lemma \ref{arrows},
$B \overline {\bf h}^*$ and $\bar A \overline {\bf h}^*$, and
\begin{equation}\label{antip} \overline{( A {\bf h})}^*= B
\overline {\bf h}^*, \quad \quad \overline{ \bar B {\bf h}^*}=\bar
A \overline {\bf h^*}.
\end{equation} If we consider the successive points of the cycle,
after ${\bf h}$ and  $A {\bf h}$, we must meet the point
$\overline {\bf h}^*$. To every point ${\bf f}= T {\bf h}$,
successive to $ {\bf h}$, there correspond, in the antipodal part
of the cycle from ${\bf h}$ to $\overline {\bf h}$, a symmetric
point, $\overline{\bf f}= \check T \overline {\bf h}$, where
$\check T$ is obtained by $T$ replacing each $A$ by $B$, and each
$B$ by $A$. So, we write finally
\[ \overline {\bf h}^*= M {\bf h}, \quad \quad  {\bf h}=  \check M \overline {\bf h^*}.  \]
Note that the last operator of $M$ must be the same as the initial
(in this case  $A$), because of eq. (\ref{antip}).

So we obtain
\[   {\bf h} = \check M  M  {\bf h}, \]
where $M$, because of the preceding remark, has an odd number $q$
of groups of operators $A$ and $B$. The operator $T=\check M M$
defining the cycle is  thus
\[  T= B^{a_1}A^{a_2}\cdots B^{a_q}A^{a_1}B^{a_2}\cdots A^{a_q}.  \]
Hence the roots $\xi^\pm( {\bf h})$ have  odd periods, non
palindromic,
\[ [a_1,a_2,\dots, a_q], \quad [a_q,a_{q-1},\dots, a_1]. \]
\hfill $\square$

{\it Remark.} In this case  $\Pi=\Gamma^2$ and $p=2P=2q$.

\subsubsection*{Proof of Theorem \ref{tII-7}} (See Figure
\ref{symm}-V).

The asymmetric case is the simplest. The number $p$ of groups of
operators  $A$ and $B$, factors of the  operator $T\in \mathcal
T^+$ satisfying $T {\bf h} = {\bf h}$ (${\bf h}$ being a principal
point of the cycle), is necessarily even.   Since there are no
symmetries, the periods of the roots $\xi^{\pm}( {\bf h})$ are non
symmetric and contain thus $p$ (even) elements.   \hfill $\square$

{\it Remark.} In this case  $\Pi=\Gamma$ and $p=P$.

\subsubsection*{Proof of Corollary \ref{corII-2}} By Theorem \ref{tII-2bis},
the sequence $s$   defines uniquely a primitive class iff this
class is invariant under antipodal symmetry.  Such a class is
either supersymmetric or antisymmetric. By Theorems
\ref{tII-3}--\ref{tII-7}, these cases are the only cases in which
the periods are odd. \hfill $\square$

\subsubsection*{Proof of Theorem \ref{nopol}} Suppose that the cycle
in $H^0$ have the symmetry of a regular n-gon, with ${\rm n}>2$.
This means that there exists  an operator $M$ of $\mathcal T^+$
and a  point ${\bf h}\in H^0$ such that ${\bf h}=M^{\rm n}{\bf
h}$.     By the results of Part I and Theorem \ref{tII-1}, all
points ${\bf h}_i$ of the cycle satisfy $\widetilde {M^{\rm
n}}_i{\bf h}_i={\bf h}_i$, for some $\widetilde {M^{\rm n}}_i$,
obtained  from   $M^{\rm n}$ by a cyclic permutation of its
factors. Among them, there is a principal point $\bf f$ such that
$\xi^+({\bf f})$ is immediately periodic
\[ \xi^+({\bf f})=[[(a_1, \dots, a_p)_1,(a_1\dots
a_p)_2,\dots,(a_1,\dots,a_p)_{\rm n}]],
\]
 being  for it $\widetilde M^{\rm
n}=(B^{a_1}A^{a_2}\dots A^{a_p})^{\bf n}$, with $p$ even. But
$\xi^+({\bf f})$ satisfies
\[ \mu(\xi^+({\bf f}))=\xi^+({\bf f}),\] $\mu$ being obtained from
$(B^{a_1}A^{a_2}\dots A^{a_p})$ by translating $A$ into $\alpha$
and $B$ into $\beta$. By Lemmas \ref{corrisp} and \ref{lemab},
also $\bf f$ satisfies
\[ {\bf f}=(B^{a_1}A^{a_2}\dots A^{a_p}){\bf f},   \]
and thus belongs to a cycle of $p$ elements.  Since a point of the
cycle cannot be visited twice,  no other points   can belong to
the same cycle, so that the cycle has no the symmetry of a regular
n-gon. \hfill $\square$

\subsubsection*{Proof of Theorem \ref{tII-8}.}
By Theorem 4.13 of Part I, the number of points inside $G_A$ and
$G_{\bar A}$  of class $C(m,n,k)$ is equal to the number $t_B$ of
factors of type $B$ in the operator $T\in \mathcal T^+$, which
defines the cycle in $H^0$.  Theorem \ref{tII-1} establishes a
correspondence between  $\Pi=[a_1,a_2,\dots,a_p]$   and the
sequence of factors $A$ and $B$ of $T$. Hence $t_B$ is equal
either to $t_{odd}$ or to $t_{even}$, depending whether  the
generators to the powers $a_j$, with $j$ odd, are $B$ or $A$.

If the period $\Gamma$ is even and palindromic, $P$ is even and
$a_i$ is equal to $a_{P+1-i}$.   Hence, when $i$ is even,
$(P+1-i)$ is odd, and vice versa. Hence  values $t_{odd}$ and
$t_{even}$, given by eq. (\ref{sum}), coincide, and being their
sum equal to $t$,  are equal to $t/2$.

If the period is even and non palindromic, $p=P$, and  values
$t_{odd}$ and $t_{even}$, given by eq. (\ref{sum}), do not
coincide necessarily.

If the period $\Gamma$ is odd, $P$ is odd. In $\Pi$, by
definition, $a_{i+P}=a_i$, for all $i=1,\dots ,P$ and when $i$ is
even, $(i+P)$ is odd and vice versa. Also in this case
$t_{odd}=t_{even}$.

On the other hand, the equality $t_{odd}=t_{even}$ do not imply,
evidently, the symmetry of the period.

\hfill $\square$

\subsection{Classes representing 1}

\subsubsection*{Proof of Corollary \ref{corII-1}}
We prove this corollary\footnote{Alternative proof: The  class of quadratic form $C(1,n,k)$ is the sole class, with
discriminant $k^2-4n$, representing 1. This class  is the group
identity of the class group. Since the inversion in the class
group corresponds to the conjugation, the identity class is
self-conjugate, and hence it is invariant under reflection with
respect to the axis $k=0$. The sole classes having this symmetry
are the supersymmetric classes (which have period odd and
palindromic) and the $k$-symmetric classes (which have period even
and bipalindromic).}   by showing that  the class of  forms representing 1
is either $k$-symmetric or supersymmetric.  The  statements will follow from Theorems 
\ref{tII-5} and \ref{tII-4}.  Indeed, in such a class  there is a representant
$\bf f$ with coefficients $(1,n,k)$, i.e., a form 
\[  f=x^2+kyx+ny^2.\]
Let  $T=\begin{pmatrix}1  & -k \\ 0 & 1 \end{pmatrix} \in \SL(2,\Z)$.
If $v=(x,y)$,  $g(v)=f(T v)=x^2-kxy+ny^2$.
We observe that ${\bf g}=(1,n,-k)= \overline {\bf f}$
and hence  $\overline {\bf f}\in C(1,n,k)$.  Therefore the   
class $C(1,n,k)$  is  either  $k$-symmetric or  supersymmetric.

 \hfill $\square$

\subsection{Classes representing 0}

Classes of forms representing 0 have discriminant equal to a
square number. The roots of the corresponding quadratic equations
are rational, and have thus a finite continued fraction.

\subsubsection*{Proof of Theorem \ref{tII-9}.}
The rational number $k/m$ is equal to $\xi^+({\bf h})$, ${\bf
h}=(m,0,-k)$.  The quadratic form ${\bf h}$, for which $\xi^+({\bf
h})>0$ and $\xi^-({\bf h})=0$,  belongs to set $F_{\bar A}$ (see
Section 4.4  in Part I).     By Theorem 4.14 of Part I, ${\bf h}$
is one of the $k$ forms with discriminant $\Delta=k^2$ and it is
the sole form in $F_A$ of its class.

Class $C(m,0,k)$ is the conjugate of  class $C(m,0,-k)$, being
$\overline{\bf h}=(m,0,k)$.    Conjugate classes, being symmetric
with respect to the plane $k=0$,  have the same number of points
in each domain.

The form ${\bf f}=A{\bf h}$  is inside $H^0$, and satisfy $\bar A
{\bf f}\in F_{\bar A}$.   By Theorem 4.19 of Part I, it is the
starting point of a chain, containing all  $t$ points in $H^0$ of
its class. The final point of that chain represents form ${\bf g}$
satisfying: ${\bf g}=T{\bf f}$ for some operator $T\in \mathcal
T^+$, product of $t-1$ generators.  Moreover, ${\bf p}=A{\bf g}\in
F_A$, and $\xi^+({\bf p})=0$.

By hypothesis  $[a_1,a_2,\dots,a_L]$ is the even continued
fraction of $\xi^+(\bf h)$. Using Lemmas \ref{lemab} and
\ref{shift}, we thus obtain:

$\xi^+({\bf f})=[a_1-1,a_2,\dots,a_L]$, and $\xi^+({\bf g})=[1]$.

Define $\tau$ by  $\xi^+({\bf g})=\tau (\xi^+({\bf f}))$. Since
$L$ is even
\[\tau={ \beta^{a_L-1}\circ \alpha^{a_{L-1}}\cdots \beta^{a_2}\circ
\alpha^{a_1-1}}. \]

Hence operator $T$ sending  ${\bf f}$ to ${\bf g}$ is equal to
\[  T= B^{a_L-1} A^{a_{L-1}}\cdots B^{a_2}
 A^{a_1-1}.  \]

According to Theorem 4.19 of Part I, the number of points of
$C(m,n,k)$ inside every domain in $G_A$ and in $G_{\bar A}$ is
equal to the number of factors of type $B$ in $T$, and it is
therefore equal to $t_{even}$, given by eq. (\ref{eq4}), whereas
the number of points inside every domain in $G_B$ and in $G_{\bar
B}$ is equal to the number of factors of type $A$, and it is thus
equal to $t_{odd}$.

When the class is supersymmetric or $(m+n)$-symmetric, it is
symmetric with respect to the horizontal axis $(m+n)=0$.  The
number of points inside the  domains in  $G_A$ and $G_{\bar A}$
coincides with that inside  the domains in $G_B$ and $G_{\bar B}$.

Since  $t_{odd}+t_{even}=t-1$, we have
\begin{equation} \label{tminus} t_{odd}=t_{even}=(t-1)/2. \end{equation}

The number of points inside $H^0$ is thus odd, when $t_A=t_B$.
  The theorem is proved \hfill
$\square$

{\it Remark.} Eq. (\ref{tminus}) says  that, when the class is
supersymmetric or $(m+n)$-symmetric, the number of points inside
$H^0$ is odd.

\subsubsection*{Proof of Theorem \ref{t10}}

For $m=1,\dots,k-1$,  the rational number $k/m$  is the non zero
root of the equation $m\xi^2-k\xi=$, i.e., $\xi^+({\bf h})=k/m$
and  ${\bf h}=(m,0,-k)$ is the representative of the class
$C(m,0,-k)$.

Form $(m,0,k)$ is the conjugate of ${\bf h}$, and its class has
the same type of symmetry as the class of ${\bf h}$.

i)  Lemma 4.17 of Part I  proves that $C(0,0,k)$ is
supersymmetric. If $m=-k/2$, let ${\bf h}=(m,0,-2m)$. Then form
${\bf f}:=A{\bf h}=(m,-m,0)$  is the central point of $H^0$. By
Lemma \ref{arrows},   $B{\bf f}={\bf h}^*$, $A{\bf
f}=\overline{\bf h}$, $\bar B{\bf f}=\overline{\bf h}^*$, i.e.,
the four forms obtained applying the generators $A$ and $B$ and
their inverse to ${\bf f}$ are symmetric and belong to the
boundary of $H^0$. Hence ${\bf f}$ is the only point of the orbit
inside $H^0$ and the orbit is supersymmetric.  The fact that if an
orbit $C(m,0,k)$ is supersymmetric then  $m=0$ or $m=k/2$ is
proved by the following reasoning.  The initial point ${\bf i}$ of
the chain in $H^0$ is joined by two inverse arrows   to $\bar
A{\bf i}$ and $\bar B{\bf i}$, and the final point ${\bf f}$ by
two arrows   to $A{\bf f}$ and $B{\bf f }$. Since the orbit is
supersymmetric, $\bar A {\bf i }$ and $\bar B {\bf i}$ are
antipodal as well as the arrows from them, and  $A{\bf f}$ and
$B{\bf f}$ are antipodal, as well as the arrow to them. Therefore,
the only possibility in order the chain, with its four points on
the boundary of $H^0$, have the supersymmetry, is that the initial
point and the final point either both coincide with the central
point of $H^0$, or do not exist. The first case is that of
$C(k/2,0,k)$, whereas the second case is that of $C(0,0,k)$.

ii) Let $ {\bf h}=(m,0,-k)$,  $m>0$ and $m\not=k/2$.  Let,
moreover, $m/k=[a_1,\dots, a_L]$, with $L$ even.  Suppose
$C(m,0,-k)$ be $(m+n)$-symmetric. By the same arguments proving
Theorem \ref{tII-9}, we obtain
\[   \beta^{a_L}\circ \alpha^{a_{L-1}}\cdots \beta^{a_2}\circ\alpha^{a_1}(\xi^+({\bf h}))=0,  \]
whereas, by the same arguments proving Theorem  \ref{tII-3},
\[   B^{a_L} A^{a_{L-1}}\cdots B^{a_2} A ^{a_1}{\bf h}= {\bf h}^*.  \]
The final point ${\bf g}$ of the chain staring at ${\bf f}=A{\bf
h}$ must be sent by $\bar B$ to ${\bf h}^*$, because of the
symmetry of the orbit.  Every point of the chain obtained as $T
{\bf h}$,  $T\in \mathcal T^+$, must have its symmetric, as well
as every arrow, according to Lemma \ref{arrows}. In this way we
obtain a palindromic sequence of groups of operators $A$ and $B$,
which must be even, since it starts by $A$ and ends by $B$. Like
in the proof of Theorem \ref{tII-3}, we conclude that  every
$(m+n)$-symmetric chain contains one (only one in this case)
selfadjoint point $(m+n=0)$ in $H^0$. On the other hand, if we
have  $\xi^+(m/k)$ even and palindromic, using Lemma \ref{shift}
we reach zero, i.e., a point of the boundary of $H^0$, and by the
symmetry of the corresponding operator of $\mathcal T^+$, we
conclude that this point is ${\bf h}^*$, and hence the orbit is at
least $(m+n)$-symmetric.  Such orbit could be, in fact,
supersymmetric, but this is excluded by point (i).

 iii) Let $ {\bf h}=(m,0,-k)$,
$m>0$ and $m\not=k/2$.  Let, moreover, $m/k=[a_1,\dots, a_L]$,
with $L$ odd.  Suppose $C(m,0,-k)$ be $k$-symmetric. By the same
arguments proving Theorem \ref{tII-9}, we obtain
\[   \alpha^{a_L}\circ \beta^{a_{L-1}}\cdots \beta^{a_2}\circ\alpha^{a_1}\xi^+=0,  \]
whereas, by the same arguments proving Theorem  \ref{tII-4},
\[   A^{a_L} B^{a_{L-1}}\cdots B^{a_2} A ^{a_1}{\bf h}= \overline{\bf h}.  \]
Indeed, the final point ${\bf g}$ of the chain starting at ${\bf
f}=A{\bf h}$ must be sent by $A$ to $\overline{\bf h}$, because of
the symmetry of the orbit.  Every point of the chain obtained as
$T {\bf h}$,  $T\in \mathcal T^+$, must have its symmetric, as
well as every arrow, according to Lemma \ref{arrows}. In this way
we obtain a palindromic sequence of groups of operators $A$ and
$B$, which must be odd, since it starts by $A$ and ends by $A$. On
the other hand, if  $\xi^+(m/k)$ is odd and palindromic, using
Lemma \ref{shift} we reach zero, i.e., a point of the boundary of
$H^0$, and by the symmetry of the corresponding operator of
$\mathcal T^+$, we conclude that this point is $\overline{\bf h}$,
and hence the orbit is at least $k$-symmetric. Such orbit could
be, in fact, supersymmetric, but we must exclude this case because
of point (i).

iv) For an antisymmetric orbit we should have $t_A=t_B$, because
of the complementarity between $G_A$ and $G_B$, as well as between
$G_{\bar A}$ and $G_{\bar B}$.  By Theorem \ref{tII-9}, the number
of points in $H^0$ should thus be odd.  An antisymmetric chain in
$H^0$ containing an odd number of points must contain the centre
of $H^0$, but in this case the orbit is supersymmetric, by point
(i).

v) The remaining possibility is that the chain as well as the
orbit $C(m,0,k)$ is asymmetric.  Hence this happens  iff neither
the odd nor the even continued fraction of $k/m$  are palindromic.

\hfill $\square$

\section{Reduction theory}\label{reduc}
Theorem \ref{tII-1}  has as  natural consequence a reduction
procedure (probably already known), described by the following
corollary.  We will see moreover in this section how this
reduction method is related to another  method, to which we refer
as 'classical reduction theory' (\cite{sar1}, \cite{sar2},
\cite{zag} ).

Given any indefinite form $f=mx^2+ny^2+kxy$, for which $mn\ge0$,
the procedure to  transform  it in a form of the same class
$m'x^2+n'y^2+k'xy$ such that $m'n'\le 0$ is given by the following
corollary.

\begin{cor}
Let ${\bf f}=(m,n,k)$, with $mn\ge 0$ and  $\tilde{\bf f}:=(\tilde
m,\tilde n,\tilde k)$ be obtained applying to ${\bf f}$ one of the
involutions, such that $\xi^+(\tilde {\bf f})>0$.  Let
$[\alpha_0,\alpha_1,\dots,\alpha_N,[a_1,\dots, a_p]]$ represent
the continued fraction of the  root $\xi^+(\tilde {\bf f})$ (the
possible absence of the period indicates that the root is
rational).  Then the form
\[  {\bf f'}:=(m',n',k') = \widetilde {L({\bf \tilde f}} )\]
satisfies  $m'n'\le 0$, where $L= C^{\alpha_N}\cdots A^{\alpha_2}
\ B^{\alpha_1} A^{\alpha_0}$,  being $C=A$ if $N$ is even and
$C=B$ if $N$ is odd.
\end{cor}

{\it Proof.}  Note that the involutions do not change the product
$mn$; to find the required involution,  look at  Table 1, and
choose it depending on the  value   of $\xi^+(\bf f)$. Form
$\tilde{\bf f}$ may not belong to the class of   ${\bf f}$. For
this reason at the end of the procedure we  apply the same
involution to the form obtained  as $L\tilde{\bf f}$. Indeed,  the
orbit of two forms related by an involution, either coincide (if
the class is invariant under that involution), or have no common
elements, each element of one class being related by that
involution to an element of the other class.

 The expression of $L$ follows from Lemma \ref{shift}.  We obtain in this way, if $N$ is even, that
 \[ \xi^+(L \tilde {\bf f})=[[a_1,\dots,a_P]], \]
 otherwise
\[ \xi^+(L \tilde {\bf f})= [0,[a_1,\dots,a_P]]. \]

The condition $\xi^+(\tilde {\bf f})>0$ means that $\tilde{\bf f}$
is either in $G_{\bar A}$ or in  $G_{\bar B}$  (where $m>0$ and
$n>0$). The operator $L$ belongs indeed to $\mathcal T^+$, and its
sequence of factors of type $A$ and $B$ can be read as a path from
$\tilde{\bf f}$ towards $H^0$, lying entirely  in $G_{\bar A}$ or
in  $G_{\bar B}$.   This path is unique because of Theorem 4.2 of
Part I.

Being $\xi^+(L \tilde {\bf f})$ immediately periodic, we conclude
that  $L \tilde {\bf f}\in H^0$ ($\xi^+(L\tilde {\bf f})=0$ if the
period is absent) and hence that the  form $L\tilde {\bf f}$
satisfies  $\tilde m' \tilde n'<0$ (or $\tilde m' \tilde n'=0$).
Again, the involution will not change the product $mn$, and we
obtain that form  ${\bf f'} = \widetilde {L \tilde {\bf f}}$
belongs to the same class of ${\bf f}$. \hfill $\square$

{\it Remark.} The operator $L'$ of $\PSL(2,Z)$ such that ${\bf f'}
= L'({\bf f})$ is deduced from $L$ using Lemma $\ref{arrows}$. It
 is still an operator of $\mathcal T^+$ or of $\mathcal T^-$,
whose factors indicate the unique  path from ${\bf f}$ towards
$H^0$ which lies entirely either in $G_{A}$ or in  $G_{B}$ or in
$G_{\bar A}$ or in $G_{\bar B}$.

The number of forms of that class satisfying $mn<0$ is given by
Theorem \ref{tII-8}.

\subsection{Classical reduction theory seen in the de Sitter
world}$\ \ \ \ \ \ \  \ $

{\bf Definition.} The sequence (finite or infinite) of integers
$b_i\ge 2$, $(b_0,b_1,b_2,\dots )$,  denotes the {\it modular
fraction} of the the number $\xi$:
\begin{equation}\label{mod1}  \xi = b_0-\frac{1}{b_1-\frac{1}{b_2-
\frac{1}{\cdots}}}. \end{equation}

An irrational number has infinite modular fraction, and periodic
modular fractions (with period of length $t$) are defined exactly
as the periodic continued fractions (see  Definition 3), and are
denoted by
\begin{equation}\label{mod3}  (b_0,b_1,\dots,b_{N-1},(b_N, b_{N+1},\dots,
b_{N+t-1})).\end{equation}

The following theorems represent a synthesis of the  classical
reduction theory. The proofs that we give here show a geometrical
description of this theory  in the  Sitter world.

{\bf Definition.}  An indefinite form ${\bf h}=(m,n,k)$ such that
$C(m,n,k)$ does not represent zero is said to be {\it reduced} iff
satisfies $m>0$, $n>0$, $k<0$, $m+n<|k|$.

\begin{thm}\label{reduct1} Let ${\bf h}=(m,n,k)$ be a reduced form.
Then the modular fraction   of $\xi^+({\bf h})$ is immediately
periodic, and the number of its elements is equal to the number of
reduced  forms of the class $C(m,n,k)$.  The first roots $\xi^+$
associated to the  other reduced forms are given by the cyclic
permutations of the elements of $\xi^+({\bf f})$.
\end{thm}

\begin{thm}\label{reduct2} Let ${\bf f}=(m,n,k)$ satisfy
$m>0$, $n>0$, $k<0$.  Then the modular fraction   of $\xi^+({\bf
f})$ is periodic:
\[    \xi^+({\bf f})=(b_0,b_1,\dots, b_{M},(c_1,c_2,\dots,c_t) )\]
and the operator
\[    L=RA^{b_{M}}\cdots RA^{b_2}RA^{b_1}RA^{b_0} \]
satisfies:  $L{\bf f}={\bf h}$, where ${\bf h}$ is reduced.
\end{thm}

{\it Proof of Theorem} \ref{reduct1}. The proof follows from the
results of this works using the following

\begin{lem}\label{mod} Let ${\bf f}=(m,n,k)$ and $\xi^+({\bf
f})=(a,b,c,d, \dots)$.  Then
\[ \xi^+(RA^{a}{\bf f})= (b, c,d,\dots ); \]
i.e., the cancellation of the first element  in the modular
fraction of $ \xi^+({\bf f})$ corresponds to the   action on ${\bf
f}$ of the operator $A$, iterated a number of times equal to that
element, followed by the action of $R$.
\end{lem}

{\it Proof.} By Lemma \ref{lemab},
\[ \xi^+(A^{a}{\bf f})=\alpha^a(\xi^+({\bf
f})=-\frac{1}{b-\frac{1}{c-\frac{1}{d-\frac{1}{\cdots}}}},
\]

\[ \xi^+(RA^{a}{\bf f})=\sigma \circ \alpha^a(\xi^+({\bf f}))=
b-\frac{1}{c-\frac{1}{d-\frac{1}{\cdots}}}.
\]
\hfill $\square$

 A reduced form ${\bf h}$ of a class $C(m,n,k)$ non representing
zero  belongs to $H_{\bar A}$. By Lemma \ref{corrisp}, it is
represented in $\Xi$ by a point with $\xi^+>1$ and $0<\xi^-<1$.

\begin{figure}[h]
\centerline{\epsfbox{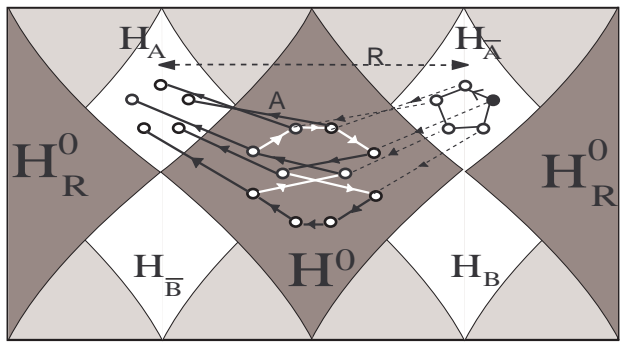}} \caption{Relation between a cycle
in $H_{\bar A}$, represented by [1,1,3,1,1,3] and a cycle in
$H^0$, represented by (3,5,3,2,2). The black dot in $H_{\bar A}$
corresponds to the form ${\bf h}$ in the proof of Theorem
\ref{reduct1} and in the Example. }\label{desit3}
\end{figure}

We know that point  $A{\bf h}$ is in $H^0$ (see Figure
\ref{desit3}). We consider now all the successive points $A^2{\bf
h}$, $A^3{\bf h}$, etc that belongs to $H^0$.  At some iteration
of order $c_1$ point $A^{c_1}{\bf h}$ exits  from $H^0$ and lies
necessarily in $H_{A}$. At this point we apply operator $R$, going
back to $H_{\bar A}$. The point ${\bf h}'=RA^{c_1}{\bf h}$
coincides with ${\bf h}$ only if ${\bf f}$ contains a sole point.
Indeed, if $RA^{c1}{\bf h}={\bf h}$, consider ${\bf f}=A{\bf h}$,
which is in $H^0$.  It satisfies $ARA^{c1}{\bf f}={\bf f}$. Using
the relation $R=\bar  A B\bar A$, we write
\[   B A^{c_1-1}{\bf f}={\bf f},\]
so concluding that the cycle in $H^0$ contains only two principal
points and, by Theorem \ref{tII-8}, that in $H_{\bar A}$ the orbit
has only one element.  If ${\bf h'}$ is different from ${\bf h}$,
we repeat the procedure till
\[   RA^{c_t}RA^{c_{t-1}}\cdots RA^{c_1} {\bf h}={\bf h}.\]
By the same reason, this  occurs  when all points in $H_{\bar A}$,
as well as in $H^0$, are visited (see figure). By Lemma \ref{mod},
using the same arguments of the proof of Theorem \ref{tII-1}, we
obtain that the modular fraction of $\xi^+(\bf h)$ is periodic,
namely
\[  \xi^+({\bf h})=((c_1,c_2,\dots, c_t)).\]

We have thus proved that the modular fraction of the first root of
an equation corresponding to a form in $H_{\bar A}$ is immediately
periodic and  the length $t$ of  its  period is equal to the
number of points of the class  in $H_{\bar A}$, i.e., the number
of reduced forms. \hfill $\square$

{\it Proof of Theorem \ref{reduct2}. } We observe now that a form
{\bf f} satisfying the conditions of the theorem  lies in $G_{\bar
A}$. Hence there is a form ${\bf h }$ in  $H_{\bar A}$ such that
${\bf f}=T{\bf h}$, where $T\in \mathcal T^{-}$. We thus write
${\bf h}=T^{-1}{\bf f}$, where $T\in \mathcal T^{+}$ is a product
of operators $A$ and $B$. We rewrite $T$ as product of operators
$A$ and $R$ by the following procedure: we replace each operator
$B$ in $T^{-1}$ by $A\bar A B \bar A A$. Hence  we replace $\bar A
B \bar A$ by $R$. We obtain a sequence containing the same factors
of type $A$ as $T^{-1}$ and
  where every factor of  type $B^i$ is
replaced by $\underbrace{(A R A )(A R  A)\cdots (A R A )}_{i \
times}$. Since the modular fraction of $\xi^+({\bf h})$ is
immediately periodic (say, $((c_1,c_2,\dots, c_t))$, by Theorem
\ref{reduct1}, using Lemma \ref{mod} we conclude that
\[T^{-1}= RA^{b_M} RA^{b_{M-1}},\cdots, RA^{b_1} RA^{b_0}  \]
iff
\[    \xi^+({\bf f})=(b_0,b_1,\dots, b_{M},(c_1,c_2,\dots,c_t) ).\]
\hfill $\square$

We obtain the following corollary:

\begin{cor}  If the period $\Pi$ of the continued
fraction of $\xi^+({\bf h})$, ${\bf h}\in H_{\bar A}$   is
$\Pi=[a_1,a_2,\dots,a_p]$, then the period of  the modular
fraction of  $\xi^+({\bf h})$ is obtained in the following way:
\begin{enumerate}
\item replace every  element $a_{2i+1}$ of the continued
fraction  having {\sl odd} index by $a_{2i+1}+2$;
\item replace every element $a_{2i}$ of the continued
fraction having {\sl even} index by a sequence containing a number
of "2" equal to $a_{2i}-1$.
\end{enumerate}
\end{cor}

{\it Remark.}  The {\sl length} of period of the modular fraction
gives the number of reduced forms in $H_A$, whereas the {\sl sum}
of all elements in the period of the continued fraction gives the
number of elements in $H^0$.

{\sc Example}. The supersymmetric class represented in Figure
\ref{desit3} is that of the Example of Theorem  \ref{tII-4}. Let
${\bf f}=(2,-1,-3)$. We have $\xi^+({\bf f
})=3/4+\sqrt{17}/4=[[1,1,3]]$. Hence $\Pi=[1,1,3,1,1,3]$. There
are indeed 10 points inside $H^0$. Point ${\bf h}=\bar{A}{\bf
f}=(2,4,-7)$ is in $H_{\bar A}$  and $\xi^+({\bf h
})=7/4+\sqrt{17}/4$ has modular fraction
\[  \xi^+({\bf h})=((3,5,3,2,2)).  \]
There are indeed 5 points in $H_{\bar A}$.

 We are now able to give the
proof of Theorem \ref{t3}.

\subsubsection*{Proof of Theorem \ref{t3}}  By the  corollary above,
given the period $(c_1,c_2,\dots, c_t)$ of a modular fraction, we
associate to it the period   $[a_1,a_2,\dots, a_p]$ of the
corresponding continued fraction, replacing every element $c_i$
greater than $2$ by an element $a_j$ equal to $c_i-2$,  and
substituting the group of $r\ge 0$ successive elements
$c_{i+1},c_{i+2},\dots,c_{i+r}$ equal to $2$, by an element
$a_{j+1}$ equal to $r+1$ (hence $a_{j+1}=1$ if $c_{i+1}>2$).  We
obtain  a period of $p$ elements with $p$ even, which is in fact
the the period $\Pi$ of the continued fraction.  By Theorem
\ref{tII-8} the sum of the elements of index odd of the continued
fraction is equal to  the sum of elements of index even if the
class has the mentioned symmetries, and is also  equal to the
number of points in $H_{A}$ and in $H_{\bar A}$ of the orbit. By
Theorem \ref{reduct1} the number of  points in $H_{\bar A}$ of the
orbit equals the number $t$ of elements of the period of the
modular fraction considered.  Hence by the equation
\[ \sum_{i=1}^t (c_i-2)=t \] we obtain  the statement of the
theorem. \hfill $\square$

{\it Remark.}  For a non symmetric class we have to consider the
modular fractions of the two inverse periods, that are different.

{\bf Acknowledgements}

I am grateful to Vladlen Timorin for having told to me the
classical reduction theory  and for his interest in this work. I
apologize number theorists for my own terminology on symmetries
and continued fractions: I hope the risk of confusion  be anyway
low, the notions I use being  all elementary.

\vfill

\newpage

\section*{Appendix: Tables of classes of indefinite forms with small discriminant  }

{\small The following  tables contain  a representative of every class
with $\Delta \le 100$;  moreover,   if $\Delta$ is not a square
integer, we give  the period $\Gamma$, its length $P$, the numbers
$t^\uparrow$ of points inside each  domain in $G_A$ and $G_{\bar
A}$,   the number $t^\downarrow$ of points inside each domain in
$G_B$ and $G_{\bar B}$, and  the type of symmetry. A star
indicates that the class is non primitive, i.e., is obtained by a
primitive class by multiplication by an integer greater than 1.

If $\Delta$ is equal to a square number, then instead of $\Gamma$
and $P$ we give the continued fraction  of $k/m$, its length $L$,
the number $t$ of points in $H^0$ and in $H^0_R$.

{\it Remark.} The classes of forms non representing zero are
either supersymmetric or $k$-symmetric if $\Delta \le 100$.

Indeed, the  first class  (i.e., with minimal discriminant) non
representing zero  having the $(m+n)$-symmetry is that in the
Example of Theorem \ref{tII-3}. It has  period $[1,2,2,1]$ and
discriminant 221.  The first  antisymmetric class has period
$[1,2,3]$ and discriminant 148, and that  asymmetric has period
$[1,1,2,3]$ and discriminant 396.

The tables show that there  are classes of forms representing zero
with $\Delta\le 100$ that have  all the possible types of
symmetry.

Moreover, we plot  in the following figure the fractions of the total
number of classes with  a  given discriminant, versus the discriminant $<10^4$, corresponding to
 the different types of  symmetries (indicated  by different symbols).
}

\begin{figure}[h]
\centerline{\epsfbox{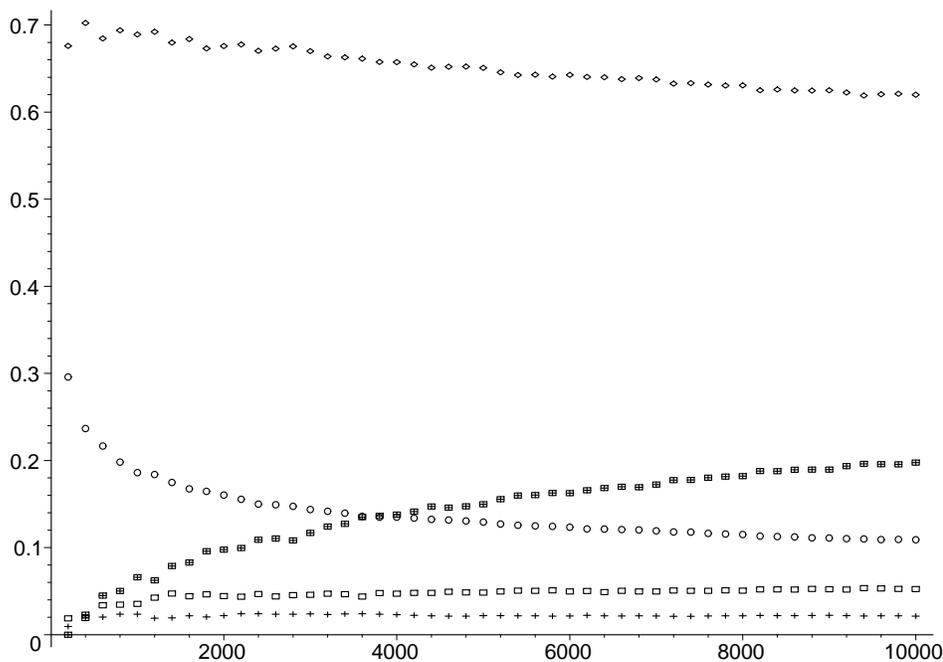}} \caption{  rhombi:  $k$-symmetric;    circles: supersymmetric;
crosses:  $(m+n)$-symmetric;   windows: asymmetric;  squares: antisymmetric}\label{numbers}
\end{figure}

\newpage
{\bf  Tables of classes non representing zero  with $\Delta< 100$}
\vskip 1 cm

 \begin{tabular} {|c||c|c|c||c|c|c||c|c|}
 \hline
$\Delta$     &  $m$ &$n$ &$k$&  $\Gamma$ & $P$ &
$t^\uparrow-t_\downarrow$ & symm. & n.p.
\\ \hline
5& 1& -1& 1& [1]& 1& 1-1 & $super$  &     \\ \hline

8& 1& -1& 2& [2]& 1& 2-2  & $super$  &     \\ \hline

12& 2& -1& 2& [2, 1]& 2&  2-1  & $k$  &     \\

  & 1& -2& 2& [1, 2]& 2& 1-2 & $k$  &     \\ \hline

13& 1& -1& 3& [3]& 1& 3-3  & $super$  &     \\ \hline

17& 2& -2& 1& [1, 3, 1]& 3& 5-5 & $super$  &     \\ \hline

20& 1& -1& 4& [4]& 1& 4-4 &  $super$  &     \\

 & 2& -2& 2& [1]& 1& 1-1 &$super$ &  $\star$ \\ \hline

21& 3& -1& 3& [3, 1]& 2& 3-1 & $k$  &     \\

  & 1& -3& 3& [1, 3]& 2& 1-3 & $k$  &     \\ \hline

24& 2& -1& 4& [4, 2]& 2& 4-2 & $k$  &     \\

 & 1& -2& 4& [2, 4]& 2& 2-4 & $k$  &     \\ \hline

 28& 3& -2& 2& [1, 1, 4, 1]& 4& 5-2 &  $k$  &     \\

 & 2& -3& 2& [1, 4, 1, 1]& 4& 2-5  & $k$  &     \\ \hline

 29& 1& -1& 5& [5]& 1& 5-5 & $super$  &     \\ \hline

 32& 2& -2& 4& [2]& 1& 2-2 & $super$  &  $\star$   \\
 & 1& -4& 4& [1,4]& 2& 1-4 &     $k$  &   \\
 & 4& -1& 4& [4,1]& 1& 4-1 & $k$  &     \\ \hline

 33& 2& -1& 5& [5, 2, 1, 2]& 4& 6-4 & $k$  &     \\ \hline

& 1& -2& 5& [2, 1, 2, 5]& 4& 4-6 & $k$  &     \\  \hline

  37& 3& -3& 1& [1, 5, 1]& 3& 7-7 &  $super$  &     \\ \hline

  40& 3& -3& 2& [1, 2, 1]& 3& 4-4 & $super$  &     \\

   & 1& -1& 6& [6]& 1& 6-6  & $super$  &     \\  \hline

41& 2& -2& 5& [2, 1, 5, 1, 2]& 5& 11-11 &  $super$  &       \\
\hline

44& 2& -1& 6& [6, 3]& 2& 6-3& $k$  &      \\

 & 1& -2& 6& [3, 6]& 2& 3-6  & $k$  &     \\ \hline

45& 3& -3& 3& [1]& 1& 1-1 &$super$ &  $\star$ \\

 & 1& -5& 5& [1, 5]& 2& 1-5 & $k$  &     \\

 & 5& -1& 5& [5, 1]& 2& 5-1 &  $k$  &     \\ \hline

48& 1& -3& 6& [2, 6]& 2& 2-6 & $k$  &      \\

 & 3& -1& 6& [6, 2]& 2& 6-2 & $k$  &     \\

 & 4& -2& 4& [2, 1]& 2& 1-2 &$k$ &  $\star$ \\

 & 2& -4& 4& [1, 2]& 2& 1-2  &$k$ &  $\star$ \\ \hline

52& 3& -3& 4& [1, 1, 6, 1, 1]& 5& 10-10 & $super$  &     \\

 & 2& -2& 6& [3]& 1& 3-3   &$super$ &  $\star$  \\ \hline

53& 1& -1& 7& [7]& 1& 7-7 &   $super$  &     \\ \hline

56& 5& -2& 4& [2, 1, 6, 1]& 4& 8-2 &  $k$  &     \\

 & 2& -5& 4& [1, 6, 1, 2]& 4& 2-8& $k$  &     \\  \hline

57& 4& -3& 3& [1, 1, 3, 7, 3, 1]& 6& 7-9 & $k$  &     \\

 & 3& -4& 3& [1, 3, 7, 3, 1, 1]& 6& 9-7 & $k$  &     \\ \hline

 \end{tabular}

\newpage

 \begin{tabular}{|c||c|c|c||c|c|c||c|c|} \hline

$\Delta$     &  $m$ &$n$ &$k$&  $\Gamma$ & $P$ &
$t^\uparrow-t_\downarrow$ & symm. & n.p.
\\ \hline

60& 2& -3& 6& [2, 3]& 2& 2-3 &  $k$  &     \\

 & 3& -2& 6& [3, 2]& 2& 3-2& $k$  &     \\

    & 1& -6& 6& [1, 6]& 2& 1-6 & $k$  &     \\ \hline
   & 6& -1& 6& [6, 1]& 2& 6-1 & $k$  &     \\ \hline
  61& 3& -3& 5& [2, 7, 2]& 3& 11-11 & $super$  &     \\ \hline
     65& 4& -4& 1& [1, 7, 1]& 3& 9-9 & $super$  &     \\
     & 2& -2& 7& [3, 1, 3]& 3& 7-7& $super$  &     \\ \hline
  68& 1& -1& 8& [8]& 1& 8-8 & $super$  &     \\
  & 4& -4& 2& [1, 3, 1]& 3& 5-5 &$super$ &  $\star$ \\ \hline
  69& 5& -3& 3& [1, 1, 7, 1]& 4& 8-2 &  $k$  &     \\
  & 3& -5& 3& [1, 7, 1, 1]& 4& 2-8 & $k$  &     \\ \hline

  72& 3& -3& 6& [2]& 1& 2-2 &$super$ &  $\star$ \\
  & 1& -2& 8& [4, 8]& 2& 8-4 & $k$  &     \\
  & 2& -1& 8& [8, 4]& 2& 8-4& $k$  &     \\ \hline
  73& 4& -4& 3& [1, 2, 3, 1, 7, 1, 3, 2, 1]& 9& 21-21& $super$  &     \\ \hline
 76& 3& -1& 8& [8, 2, 1, 3, 1, 2]& 6& 10-7 & $k$  &     \\
   & 1& -3& 8& [2, 1, 3, 1, 2, 8]& 6& 7-10 &  $k$  &     \\ \hline
  77& 1& -7& 7& [1, 7]& 2& 1-7 & $k$  &     \\
   & 7& -1& 7& [7, 1]& 2& 7-1 & $k$  &    \\ \hline

  80& 4& -4& 4& [1]& 1& 1-1 &$super$ &  $\star$ \\
  & 2& -2& 8& [4]& 1& 4-4 &$super$ &  $\star$ \\
  & 1& -4& 8& [2, 8]& 2& 2-8 & $k$  &     \\
    & 4& -1& 8& [8, 2]& 2& 8-2  & $k$  &     \\ \hline
    84& 6& -2& 6& [3, 1]& 2& 3-1 &$k$ &  $\star$ \\
   & 2& -6& 6& [1, 3]& 2& 1-3  &$k$ &  $\star$ \\
   & 4& -3& 6& [2, 1, 1, 8, 1, 1]& 6& 4-10 &$k$  &      \\
    & 3& -4& 6& [1, 1, 8, 1, 1, 2]& 6& 10-4 & $k$  &     \\ \hline
  85& 3& -3& 7& [2, 1, 2]& 3& 5-5 & $super$  &     \\
    & 1& -1& 9& [9]& 1& 9-9  & $super$  &     \\  \hline
  88& 2& -3& 8& [2, 1, 8, 1, 2, 4]& 6& 12-6  & $k$  &     \\
   & 4& -2& 8& [4, 2, 1, 8, 1, 2]& 6& 6-12  & $k$  &     \\ \hline
  89& 4& -4& 5& [1, 1, 4, 9, 4, 1, 1]& 7& 21-21 & $super$  &     \\ \hline
  92& 1& -7& 8& [1, 3, 1, 8]& 4& 2-11  & $k$  &     \\
 & 7& -1& 8& [8, 1, 3, 1]& 4& 11-2 & $k$  &     \\ \hline
  93& 1& -3& 9& [3, 9]& 2& 3-9  & $k$  &     \\
   & 3& -1& 9& [9, 3]& 2& 9-3 & $k$  &     \\ \hline
  96& 5& -3& 6& [2, 1, 1, 1]& 4& 3-2 & $k$  &     \\
   & 3& -5& 6& [1, 1, 1, 2]& 4& 2-3  & $k$  &     \\
    & 2& -4& 8& [2, 4]& 2& 2-4 &$k$ &  $\star$ \\
    & 4& -2& 8& [4, 2]& 2& 4-2 &$k$ &  $\star$ \\
    & 1& -8& 8& [1, 8]& 2& 1-8 &  $k$  &     \\
   & 8& -1& 8& [8, 1]& 2& 8-1  & $k$  &     \\ \hline
  97& 2& -2& 9& [4, 1, 2, 2, 9, 2, 2, 1, 4]& 9& 27-27  & $super$  &     \\ \hline
   \end{tabular}

\vfill  \pagebreak \newpage
{\bf Tables of classes representing zero  with $\Delta\le 100$}

\vskip 1 cm

\begin{tabular}{|c||c|c|c||c|c|c|c||c|c|} \hline

$\Delta$     &  $m$ &$n$ &$k$&  ${k}/{m}$ & $L$ & $t$ &
$t^\uparrow-t_\downarrow$ & symm. & n.p.
\\ \hline
             1& 0& 0& 1& 0& 0& 0& 0-0& $super$ &  \\ \hline

             4& 0& 0& 2& 0& 0& 0& 0-0& $super$& $\star$ \\

               & 1& 0& 2& [2]& 1& 1& 0-0& $super$&  \\ \hline

             9& 0& 0& 3& 0& 0& 0& 0-0& $super$& $\star$ \\

               & 1& 0& 3& [3]& 1& 2& 1-0& $k$& \\

            & 2& 0& 3& [1, 1, 1]& 3& 2& 0-1& $k$& \\  \hline

            16& 0& 0& 4& 0& 0& 0& 0-0& $super$& $\star$ \\

              & 1& 0& 4& [4]& 1& 3& 2-0& $k$& \\

             & 2& 0& 4& [2]& 1& 1& 0-0& $super$& $\star$ \\

           & 3& 0& 4& [1, 2, 1]& 3& 3& 0-2& $k$& \\ \hline

            25& 0& 0& 5& 0& 0& 0& 0-0& $super$& $\star$ \\

              & 1& 0& 5& [5]& 1& 4& 3-0& $k$& \\

           & 2& 0& 5& [2, 2]& 2& 3& 1-1& {\it m$+$n}& \\

        & 3& 0& 5& [1, 1, 1, 1]& 4& 3& 1-1& {\it m$+$n}& \\

           & 4& 0& 5& [1, 3, 1]& 3& 4& 0-3& $k$& \\ \hline

            36& 0& 0& 6& 0& 0& 0& 0-0& $super$& $\star$ \\

              & 1& 0& 6& [6]& 1& 5& 4-0& $k$& \\

             & 2& 0& 6& [3]& 1& 2& 1-0& $k$& $\star$ \\

             & 3& 0& 6& [2]& 1& 1& 0-0& $super$& $\star$ \\

         & 4& 0& 6& [1, 1, 1]& 3& 2& 0-1& $k$& $\star$ \\

           & 5& 0& 6& [1, 4, 1]& 3& 5& 0-4& $k$& \\ \hline

            49& 0& 0& 7& 0& 0& 0& 0-0& $super$& $\star$ \\

              & 1& 0& 7& [7]& 1& 6& 5-0& $k$& \\

          & 2& 0& 7& [3, 1, 1]& 3& 4& 2-1& $asymm$& \\

          & 3& 0& 7& [2, 2, 1]& 3& 4& 1-2& $asymm$& \\

        & 4& 0& 7& [1, 1, 2, 1]& 4& 4& 2-1& $asymm$& \\

       & 5& 0& 7& [1, 2, 1, 1]& 4& 4& 1-2& $asymm$& \\

           & 6& 0& 7& [1, 5, 1]& 3& 6& 0-5& $k$& \\ \hline

            64& 0& 0& 8& 0& 0& 0& 0-0& $super$& $\star$ \\

              & 1& 0& 8& [8]& 1& 7& 6-0& $k$& \\

             & 2& 0& 8& [4]& 1& 3& 2-0& $k$& $\star$ \\

           & 3& 0& 8& [2, 1, 2]& 3& 4& 2-1& $k$& \\

           & 4& 0& 8& [2]& 1& 1& 0-0& $super$& $\star$ \\

        & 5& 0& 8& [1, 1, 1, 1, 1]& 5& 4& 1-2& $k$& \\

        & 6& 0& 8& [1, 2, 1]& 3& 3& 0-2& $k$& $\star$ \\

        & 7& 0& 8& [1, 6, 1]& 3& 7& 0-6& $k$& \\ \hline

\end{tabular}

 \begin{tabular}{|c||c|c|c||c|c|c|c||c|c|} \hline

$\Delta$     &  $m$ &$n$ &$k$&  ${k}/{m}$ & $L$ & $t$ &
$t^\uparrow-t_\downarrow$ & symm. & n.p.
\\ \hline

            81& 0& 0& 9& 0& 0& 0& 0-0& $super$& $\star$ \\

              & 1& 0& 9& [9]& 1& 8& 7-0& $k$& \\

          & 2& 0& 9& [4, 1, 1]& 3& 5& 3-1& $asymm$& \\

             & 3& 0& 9& [3]& 1& 2& 1-0& $k$& $\star$ \\

          & 4& 0& 9& [2, 3, 1]& 3& 5& 1-3& $asymm$& \\

        & 5& 0& 9& [1, 1, 3, 1]& 4& 5& 3-1& $asymm$& \\

          & 6& 0& 9& [1, 1, 1]& 3& 2& 0-1& $k$& $\star$ \\

        & 7& 0& 9& [1, 3, 1, 1]& 4& 5& 1-3& $asymm$ & \\

           & 8& 0& 9& [1, 7, 1]& 3& 8& 0-7& $k$& \\ \hline

           100& 0& 0& 10& 0& 0& 0& 0-0& $super$ & $ \star$ \\

             & 1& 0& 10& [10]& 1& 9& 8-0& $k$& \\

            & 2& 0& 10& [5]& 1& 4& 3-0& $k$& $ \star$ \\

          & 3& 0& 10& [3, 3]& 2& 5& 2-2& {\it m$+$n}& \\

         & 4& 0& 10& [2, 2]& 2& 3& 1-1& {\it m$+$n}& $ \star$ \\

            & 5& 0& 10& [2]& 1& 1& 0-0& $super$& $ \star$ \\

      & 6& 0& 10& [1, 1, 1, 1]& 4& 3& 1-1& {\it m$+$n}& $ \star$ \\

       & 7& 0& 10& [1, 2, 2, 1]& 4& 5& 2-2& {\it m$+$n}& \\

         & 8& 0& 10& [1, 3, 1]& 3& 4& 0-3& $k$& $ \star$ \\

          & 9& 0& 10& [1, 8, 1]& 3& 9& 0-8& $k$& \\ \hline
          \end{tabular}

\end{document}